\documentclass{imsart}

\RequirePackage[OT1]{fontenc}
\RequirePackage{hyperref}

\usepackage{amsfonts}
 
\usepackage{graphicx}
\usepackage{amssymb}
\usepackage{amsmath}
\usepackage{amsbsy}
\usepackage{theorem}
\usepackage{graphicx}
\usepackage{caption}
\usepackage{geometry}
\usepackage{relsize} 
\usepackage{textcomp }
 \usepackage{latexsym}

 \usepackage{diagbox} 
  \usepackage{pgfplots}
  \usepgfplotslibrary{fillbetween}
 
\usepackage{amsfonts}

\makeatletter
\newcommand*{\rom}[1]{\expandafter\@slowromancap\romannumeral #1@}
\makeatother

\makeatletter
\def\Let@{\def\\{\notag\math@cr}}
\makeatother
\usepackage[utf8]{inputenc}
\usepackage[english]{babel}
 
\theoremstyle{plain}

\startlocaldefs
\newtheorem{proposition}{Proposition}[section]
\newtheorem{Definition}{Definition}[section]
\newtheorem{assumption}{Assumption}[section]
\theoremstyle{plain}
\newtheorem{theorem}{Theorem}[section]
\newtheorem{remark}{Remark}[section]
\newtheorem{lemma}{Lemma}[section]

\newtheorem{example}{Example}[section]

\numberwithin{equation}{section}
\theoremstyle{plain}
\endlocaldefs 
 \definecolor{bisque}{rgb}{1.0, 0.89, 0.77}
\begin{document}
 
\begin{frontmatter}
\title{Estimation of seasonal long-memory parameters}
\runtitle{Estimation of seasonal long-memory parameters}
%
\begin{aug}
\author{\fnms{Huda Mohammed} \snm{Alomari}\ead[label=e2]{halomari@students.latrobe.edu.au}}
\address{Department of Mathematics and
Statistics, La Trobe University, Melbourne, 3086, Australia\\ and
Department of Mathematics, Al-Baha University, Saudi Arabia\\\printead{e2}}
\author{\fnms{Antoine} \snm{Ayache}\ead[label=e3]{Antoine.Ayache@univ-lille.fr}}\address{Laboratoire Paul-Painlev\'e (UMR CNRS 8524), Universit\'e de Lille, B\^atiment M2, Cit\'e Scientifique, 59655 Villeneuve d'Ascq, France \\\printead{e3}}  
\author{\fnms{Myriam} \snm{Fradon} \ead[label=e4]{Myriam.Fradon@univ-lille.fr}}\address{Laboratoire Paul-Painlev\'e (UMR CNRS 8524), Universit\'e de Lille, B\^atiment M2, Cit\'e Scientifique, 59655 Villeneuve d'Ascq, France \\\printead{e4}} 

\author{\fnms{Andriy} \snm{Olenko}\ead[label=e1]{ a.olenko@latrobe.edu.au}}
\address{ Department of Mathematics and Statistics, La Trobe University, Melbourne, 3086, Australia\\\printead{e1}}
 \runauthor{H. M. Alomari et al.}
\end{aug}

\begin{abstract}
This paper studies seasonal long-memory processes with Gegenbauer-type spectral densities. Estimates for singularity location and long-memory parameters based on general filter transforms are proposed.  It is proved that the estimates are almost surely convergent to the true values of parameters. Solutions of the estimation equations are studied and adjusted statistics are proposed. Numerical results are presented to confirm the theoretical findings.

\end{abstract}

\begin{keyword}
\kwd{Gaussian stochastic process}
\kwd{seasonal/cyclical long memory}
\kwd{wavelet transformation} 
\kwd{filter }
\kwd{Gegenbauer-type spectral densities}
\kwd{estimators of parameters.}
 
\end{keyword}
\tableofcontents
\end{frontmatter}

\section{Introduction}

The importance of long memory can be seen in various applications, for 
instance in finance, internet modeling, hydrology, linguistics, DNA 
sequencing and other areas, see \cite{cont:2005}, \cite{Dette:2017}, \cite{Leonenko:2013}, \cite{Parka:2011}, \cite{Piptaq: 2017}, \cite{Samo:2007}, \cite{Samo:2016}, \cite{Willinger:2003} and the references therein.

Usually, for a stationary finite-variance random process $X(t), t\in \mathbb R,$  
long memory or long-range dependence is defined as non-integrability of
 its covariance function
 $\mathsf B(r)={\rm cov}(X(t+r),X(t)),$ i.e. $\label{long}\int_{0}^\infty |\mathsf B(r)|\,dr=+\infty,$ or, more precisely, as a  hyperbolic asymptotic behaviour of $ \mathsf B(\cdot).$

It is known that the phenomenon of long-range dependence is related to
 singularities of spectral densities, see \cite{LeonOle:2013}. The majority of publications study the case when spectral densities are unbounded at the origin.
However, singularities at non-zero frequencies
 play an important role in investigating cyclic or seasonal behavior of
 time series. Two classical models in the literature to describe cyclic behaviours of time series  are
 
   (i) a sum of a periodic deterministic trend and a stationary random noise, 

  (ii) ARMA  model with a spectral peak outside the origin. 
 
A cyclic long-range  dependent process, that will be referred to as (iii), is an intermediate case between (i) and (ii) as it has a pole in its spectral density, see \cite{Giraitis:2001}.

\begin{figure}[h]
\noindent \hspace{-1cm}	\begin{minipage}{5.2cm}
		\includegraphics[width= 4.9cm, height= 3.5cm,trim=1.0cm 1cm 1cm 1cm, clip=true]{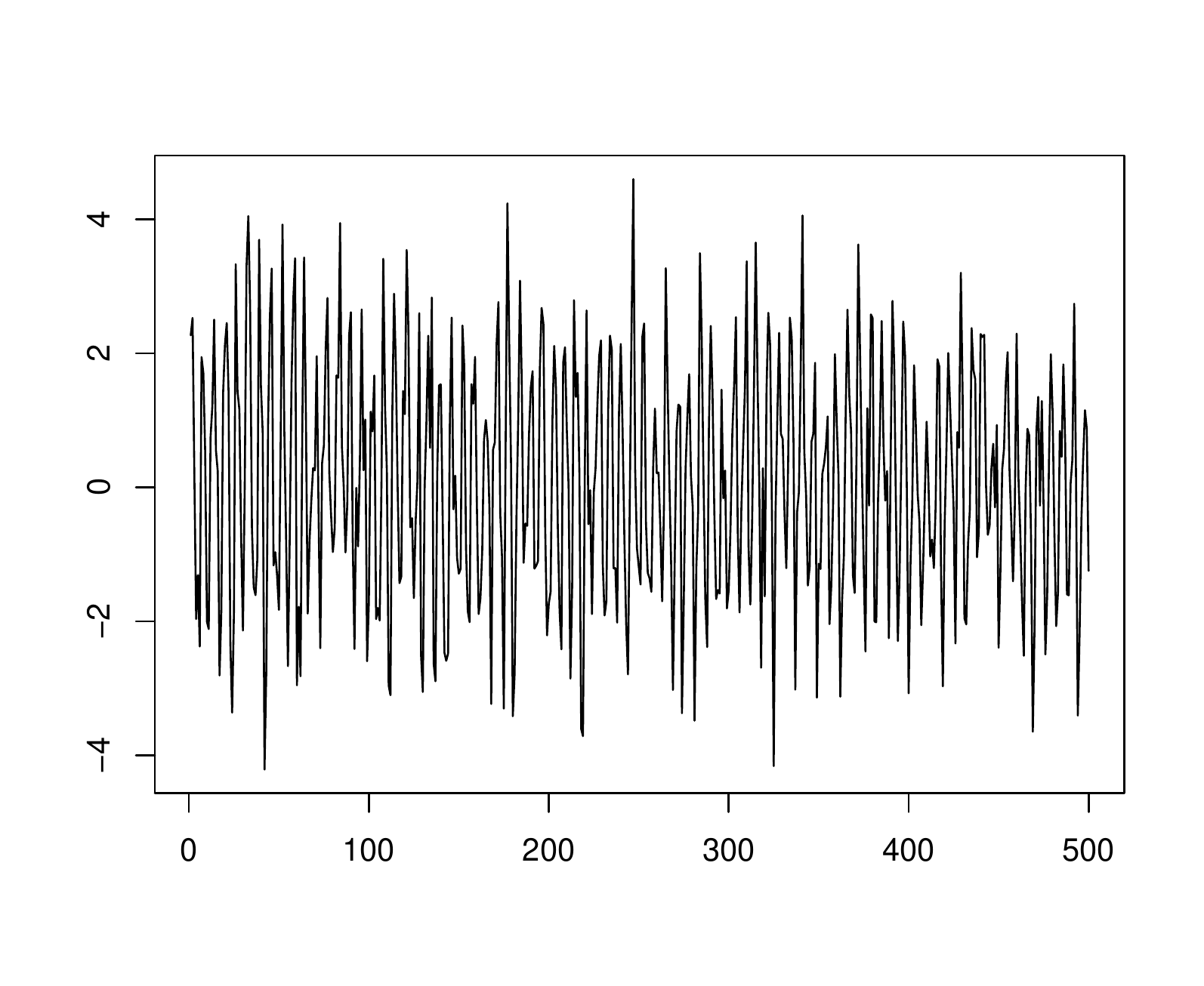}  \end{minipage}
	\begin{minipage}{5cm}
		\includegraphics[width= 4.7cm, height= 3.5cm,trim=2cm 1cm 1cm 1cm, clip=true]{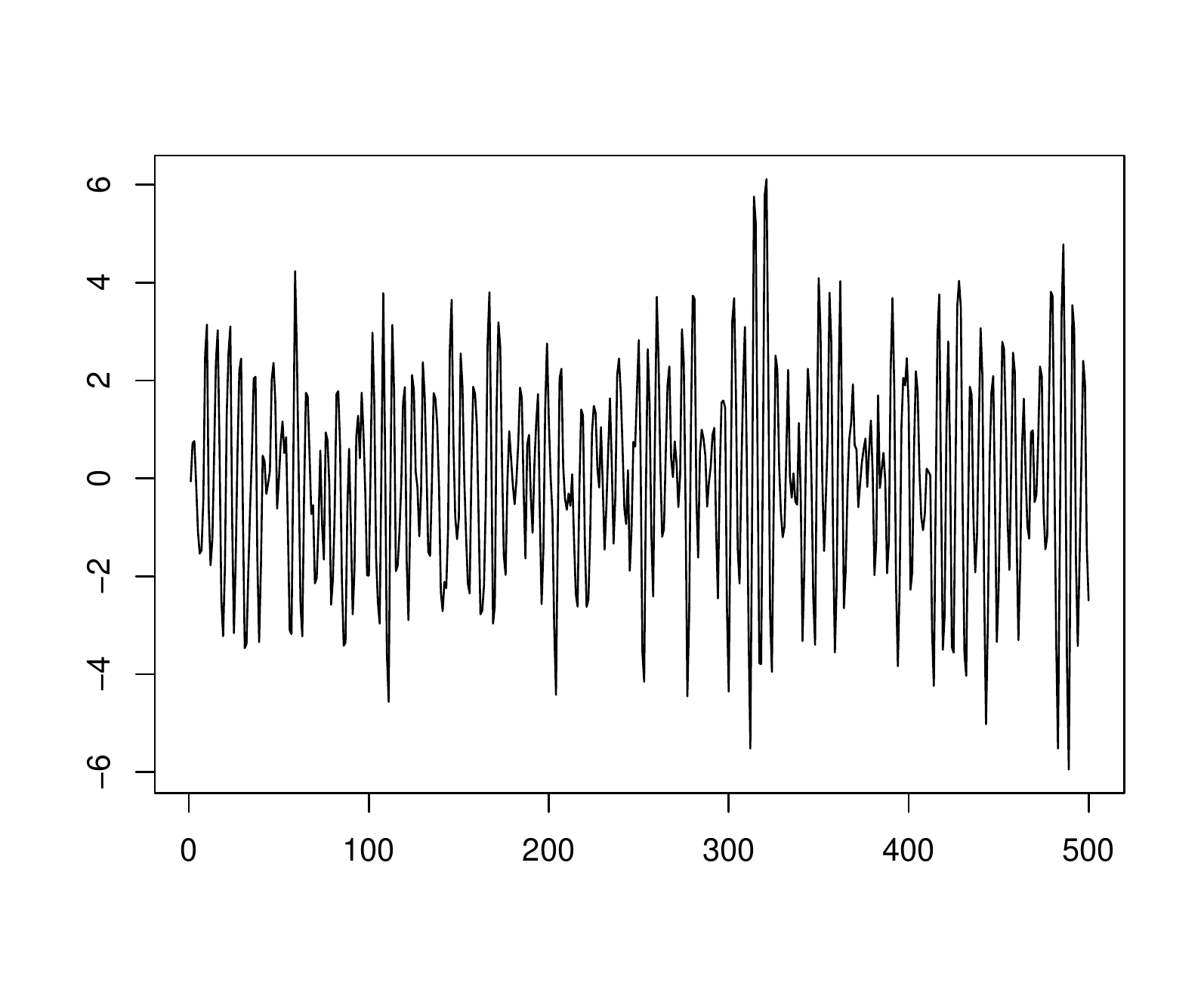}  \end{minipage}
	\begin{minipage}{4.5cm}
		\includegraphics[width= 4.7cm, height= 3.5cm,trim=2cm 1cm 1cm 1cm, clip=true]{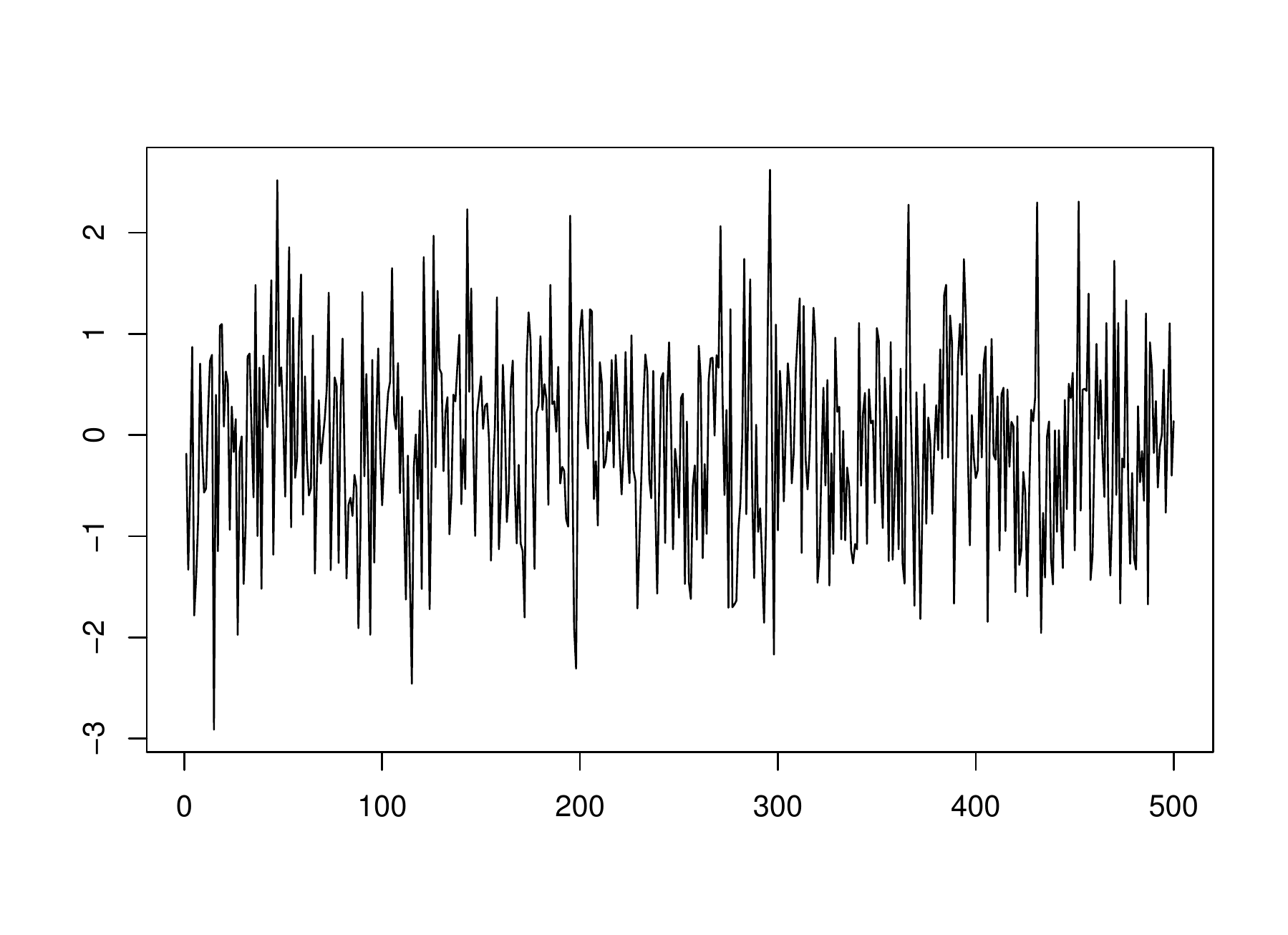}  \end{minipage}\vspace{-0.1cm}
	
\noindent	\begin{minipage}{5cm}
		\includegraphics[width=5cm, height= 4cm,trim=0.45cm 1cm 1cm 1cm, clip=true]{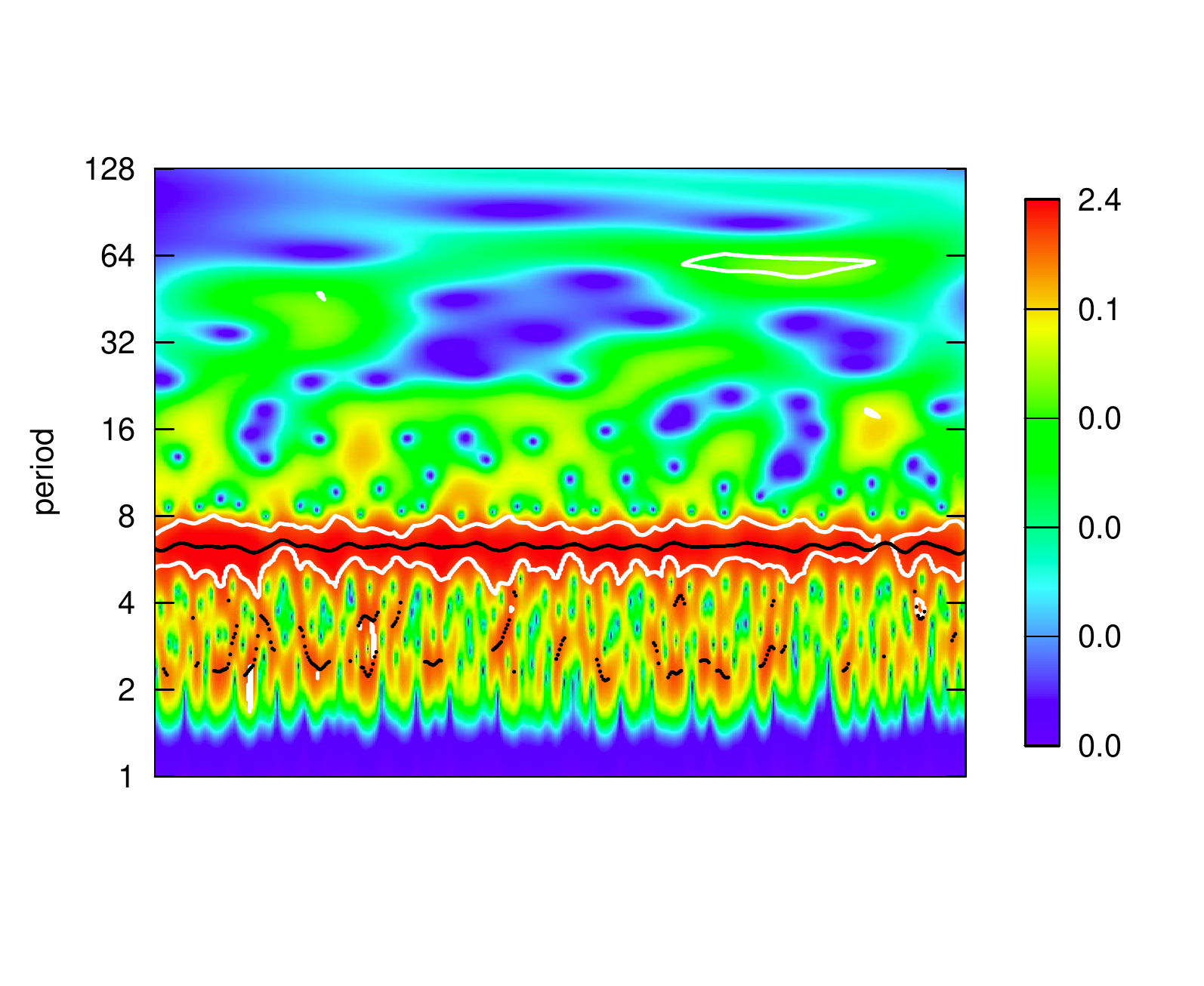}  \end{minipage}
		\begin{minipage}{5cm}
		\includegraphics[width= 4.8cm, height= 4cm,trim=2cm 1cm 1cm 1cm, clip=true]{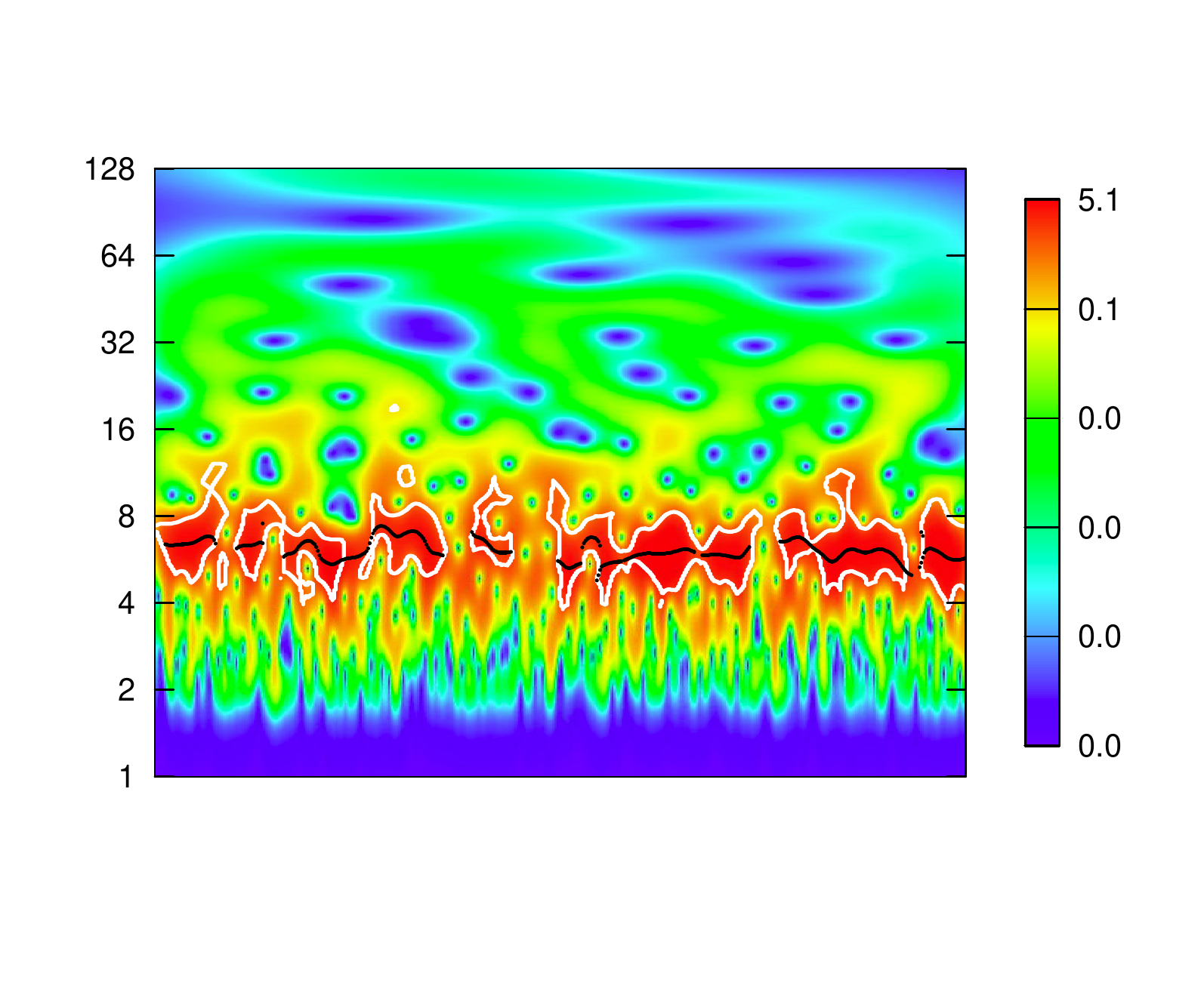}  \end{minipage}
		\begin{minipage}{4.8cm}
		\includegraphics[width=4.8cm, height= 3.4cm,trim=2cm 1cm 0.1cm 1cm, clip=true]{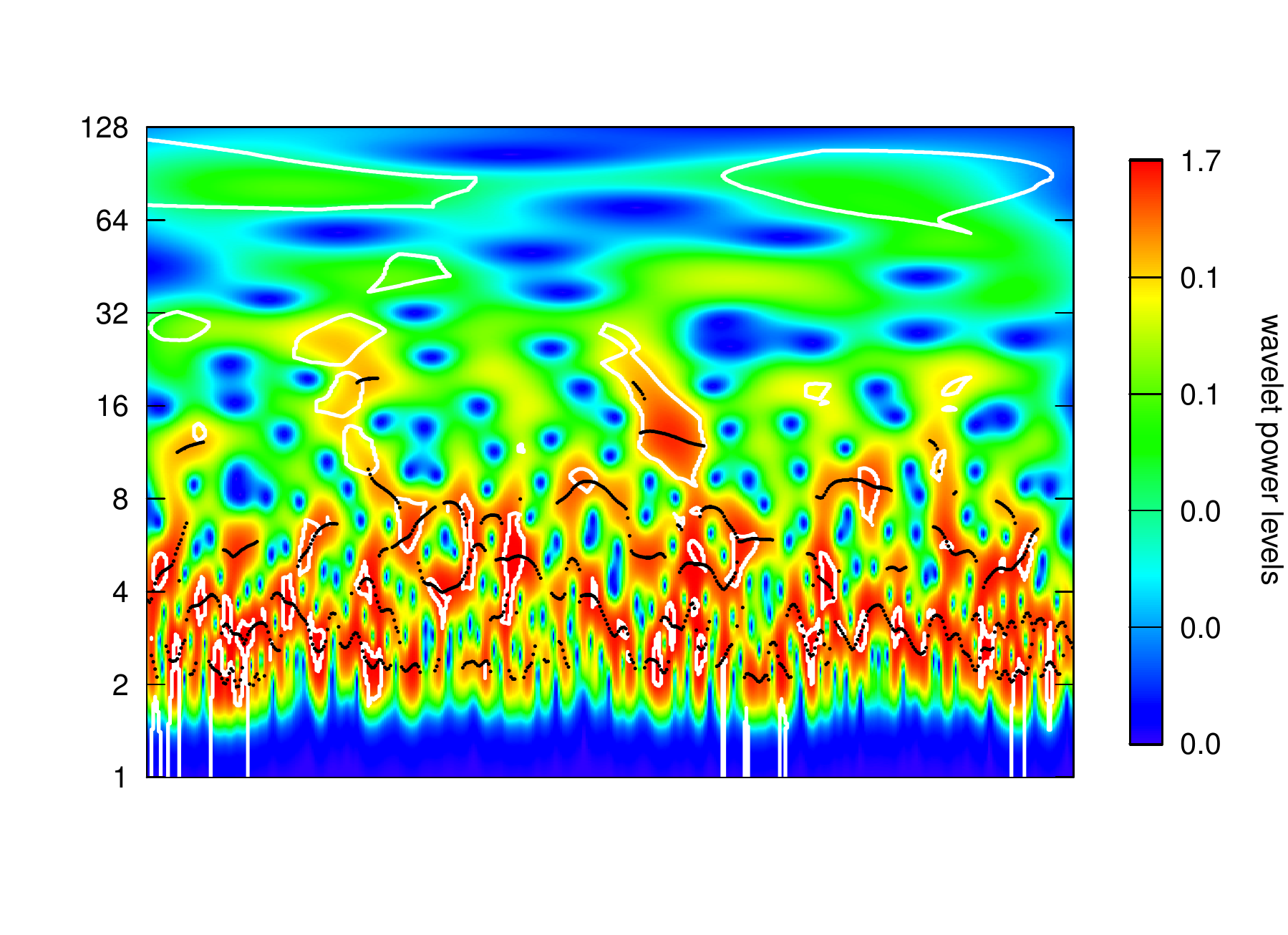}  \end{minipage}\vspace{-0.6cm}
	\caption{Three types of time series and their wavelet power spectra.}\label{fig1*}
\end{figure}

The first row of Figure~\ref{fig1*} shows realisations of models (i), (ii) and (iii) from left to right. The stochastic processes $X(t)=2\sin(t)+\varepsilon_t$ and $X(t)=0.9X(t-1)-0.8X(t-2)+\varepsilon_t,$ where $\varepsilon_t$ is a zero-mean white noise, were used as models (i) and (ii) respectively. The Gegenbauer random process from Section~\ref{sec_6} was used as model (iii)  in the simulations.  For each simulation the second row of Figure~\ref{fig1*} gives the corresponding wavelet power spectra.   It suggests that estimation of parameters might be more challenging problem for model (iii) than for cases (i) and (ii). Unexpectedly, relatively few publications on the matter are related to cyclical or seasonal long-memory processes. A survey of some recent asymptotic results for cyclical long-range dependent random processes and fields can be found in \cite{ArtRob:1999}, \cite{Ivanov:2013}, \cite{Klyolen:2012} and \cite{Olenko:2013}.  It was demonstrated in  \cite{Olenko:2013} that singularities at non-zero frequencies can play an important role in limit theorems even for the case of linear functionals.
\noindent

Several parametric and semi-parametric methods were proposed for the case when poles of spectral densities are unknown, see \cite{ArtRob:2000}, \cite{Giraitis:2001}, \cite{Hidalgo:1996}, \cite{Whitcher:2004} and the references therein.  Various problems in statistical inference of random processes and fields 
characterized by certain singular properties of  their spectral densities were 
investigated in \cite{Leonenko:1999}.  Some methods for estimating a 
singularity location were suggested by \cite{ArtRob:2000},
 \cite{Giraitis:2001} and \cite{Ferrara:2001}.
 
The asymptotic theory for Gaussian maximum likelihood estimates (MLE) of seasonal long-memory models was developed in \cite{Giraitis:2001}. The quasi-likelihood methods were studied in \cite{Hosoya:1997}.
The paper \cite{Hidalgo:1996} studied 
limit theorems for spectral density estimators and functionals with
 spectral density singularities at the origin and possibly at other 
frequencies. Some results about consistency and asymptotic 
normality of the spectral density estimator were obtained. The paper \cite{Yajima:1985} proposed the  MLE and the least squares estimator for  the long-memory time series models from \cite{GranRose:1980}
 and  the ARFIMA model in \cite{Hosking:1981}. They 
examined the consistency, the limiting distribution, and the rate of 
convergence of these estimators.  The least square estimator   
method was used in \cite{Beran:2009} to estimate the 
long-rang dependence parameter assuming that the singularity point
 is at the origin.

The minimum contrast estimator (MCE) methodology has been applied 
in a variety of statistical areas, in particular, for long-range dependent 
models.  The article \cite{Anh:2007} discussed consistency and asymptotic normality 
of a class of MCEs for random processes with short- or long-range 
dependence based on the second- and third-order cumulant spectra. In
\cite{Guo:2009} it was demonstrated that the Whittle maximum likelihood 
estimator is consistent and asymptotically normal for stationary seasonal 
autoregressive fractionally integrated moving-average
processes. Consistency and asymptotic normality of MCEs for parameters of Gegenbauer random 
processes and fields were obtain in \cite{Espejo:2015}. 
More details on the current state of  the MCE theory for long-memory
 processes and additional references can be found in \cite{Alomari:2017}. 

Unfortunately, the approaches developed in \cite{Alomari:2017}, \cite{Espejo:2015}, \cite{LeonSak:2006} for Gegenbauer-type long-memory models use specific weight functions that correspond to locations of singularities of spectral densities. Therefore,
 these methods can be applied only if locations of singularities are known
 (for example, were estimated before or determined by particular applications). 
These methods can't be applied in situations where the both long memory 
and seasonality parameters are unknown or have to be estimated simultaneously.

The article \cite{Whitcher:2004} proposed to use wavelet transforms to estimate parameters of seasonal long-memory time series. Simulation studies were used to validate the approach and to compare it with other techniques. Unfortunately, there were no rigorous studies to justify the method and establish statistical properties of the estimators, except the case of the singularity at the origin, see \cite{Clausel:2014} and the references therein. 

This research addresses  this problem and gives first steps in developing
 simultaneous estimators for the both parameters. The paper deals with  Gegenbauer-type seasonal long-memory parametric models. The  Gegenbauer spectral density $f(\cdot)$ has the following form and asymptotic  behaviour around its poles $ \pm\nu$
\begin{align*}f(\lambda)=C\left(2\left|\cos\lambda-\cos\nu\right|\right)^{-2\alpha}
= C\left(4\left| \sin\left(\frac{\lambda+\nu}{2}\right) \sin\left(\frac{\lambda-\nu}{2}\right)\right|\right)^{-2\alpha}\sim C \left|\lambda^2-\nu^2\right|^{-2\alpha}, 
\end{align*}
when $\lambda \to \pm \nu.$

The detailed review of the statistical inference theory for Gegenbauer random processes and fields can be found in \cite{Espejo:2015}. 

We use the idea from \cite{Bardet:2010} to develop the first estimation equation. Namely, we study asymptotic properties  of a filter transformation 
of seasonal long-memory processes. As a particular case this transformation
 includes wavelet transformations. To get the second estimation 
equation we propose a new approach that is based on asymptotic behaviour  
of increments of the filter transformation. Finally, we investigate properties of the solutions to the estimation equations and propose adjusted statistics for the both seasonal and long-memory parameters. 
The developed methodology includes wavelet transformations as a particular case.
Therefore, it is  potentially very useful for real applications as it can employ 
the existing wavelet methods and software, which are more powerful and 
faster than programs for numerical integration and optimization 
required by the MLE and MCE methods.
 
 The article is organized as follows. In Section~\ref{sec_2} we give basic definitions
 and notations. The first equation to estimate the 
 parameters is derived in Section~\ref{sec_3}. 
Section~\ref{sec_{4}}, further studies properties of filter transforms 
and their increments. Then these results are used to derive
 the second estimation equation. In Section~\ref{sec_{5}}, estimators of location and long memory parameters are proposed and studied.  
Simulation studies which support the
 theoretical findings are presented in Section~\ref{sec_6}.

All computations and simulations in the article were performed using the software R version~3.5.0 and Maple 17, Maplesoft.
 \section{Definitions and auxiliary results}\label{sec_2}
This section introduces classes of stochastic processes and their filter transforms that are studied in the paper.

We consider a measurable mean-square continuous stationary zero-mean 
Gaussian stochastic process $ X(t), t\in \mathbb{R},$ defined 
on a probability space $(\Omega, \mathcal{F}, P),$ with the 
covariance function
\[B(r):=Cov(X(t), X(t'))=\int_ \mathbb{R} \mathrm{e}^{iu(t-t')} F(du), \quad t, t' \in \mathbb{R}, \] 
where $r=t-t'$ and $F(\cdot)$ is a non-negative finite measure on $\mathbb{R}.$
\begin{Definition}
The random process $X(t), t \in \mathbb{R},$ is said to possess an absolutely
 continuous spectrum if there exists a non-negative function $f(\cdot)\in L_{1}( \mathbb{R})$ such that
\[F(u)=\int_{-\infty }^ {u} f(\lambda) d\lambda, \quad  u\in \mathbb{R}. \]
\end{Definition}
The function $f(\cdot)$ is called the spectral density of the process $ X(t).$

The process $ X(t), t \in \mathbb{R},$ with absolutely continuous spectrum has the 
following isonormal spectral representation
\[X(t)= \int_ {\mathbb{R}}\mathrm{e}^{i t\lambda} \sqrt{f(\lambda)} dW(\lambda),\]
where $W(\cdot)$ is a complex-valued Gaussian orthogonal random measure on $\mathbb{R}.$

For simplicity, in this paper we consider the case of real-valued $X(t).$ Therefore, 
we assume that $f(\cdot)$ is an even function and the random measure is such that 
 $W\left(\left[\lambda_1,\lambda_2\right]\right)=W\left(\left[-\lambda_2,-\lambda_1\right]\right)$ for any $\lambda_2>\lambda_1>0,$ see \S 6 in \cite{Taqqu:1979}. As all estimates 
in the paper use absolute values of integrands, the obtained results can also
be  rewritten for complex-valued processes.
\begin{assumption}\label{Assumption_1}
Let the spectral density $f(\cdot)$ of $X(t)$ admit the following representation
\[f(\lambda)=\frac{h(\lambda)}{|\lambda^{2}-s_{0}^{2}|^{2\alpha}},\quad \lambda \in \mathbb{R},\]
where $s_{0}> 0,\, \alpha\in (0,\frac{1}{2})$ and $h(\cdot)$ is an even non-negative bounded function that is four times boundedly differentiable on $[-\frac{1}{2},\frac{1}{2}].$ Its derivatives of order $i$ are denoted by $h^{(i)}(\cdot)$  and satisfy $h^{(i)}(0)=0,$ $i=1,2,3,4.$ Also, $h(0)=1,$ $h(\cdot)>0$ in some neighborhood of $\lambda=\pm  s_0,$ and for all $\epsilon>0$ it holds 
\[\int_\mathbb{R}\frac{h(\lambda)}{(1+|\lambda|)^\varepsilon }d\lambda<\infty.\]

\end{assumption}
\begin{remark}\label{remark_1}
 Stochastic processes satisfying Assumption $ \ref{Assumption_1}$ 
have seasonal long memory, as their spectral densities have 
singularities at non-zero locations $ \pm s_{0}.$ The boundedness of $h(\cdot)$ guarantees that the singularities of $f(\cdot)$ are only in $\pm s_0.$ The parameter 
 $\alpha$ is a long-memory parameter. The parameter $s_{0}$ determines seasonal or cyclic behaviour. Covariance functions of 
such processes exhibit hyperbolically decaying oscillations and $ \int\limits_ {\mathbb{R}}|B(r)|dr=+\infty$ as $\alpha\in\left(0,1/2\right),$ see \rm\cite{ArtRob:1999}.
\it{The Gegenbauer random processes have seasonal long-memory behaviour determined by Assumption~\rm\ref{Assumption_1}} \it{as their spectral densities $f(\lambda)\sim c \left|\lambda^2-s_{0}^2\right|^{-2\alpha},$ $\lambda \to \pm s_{0},$ around the Gegenbauer frequency  $s_{0},$ see \rm\cite{Chung:1966}, \rm\cite{Espejo:2015}.}
\end{remark}
\begin{remark}\label{remark_01}
 The conditions on $h(\cdot)$ guarantee that $f(\cdot)$ is a spectral density with only singularity locations at $\lambda=\pm s_0.$ The differentiability conditions on $h(\cdot)$ and its derivatives can be relaxed and replaced by H\"{o}lder assumptions in some neighborhood of the origin.
\end{remark}
The smoothness conditions guarantee the following technical inequalities required for the proof.
\begin{lemma}\label{remark_0}
 For $0\leq\lambda\leq\tilde{\lambda}\leq\frac{1}{2}$ it holds\\
\[ | h(\tilde{\lambda}) -h(\lambda)|  \leq c_1|\tilde{\lambda}^2-\lambda^2|,\]
\[h(\lambda) \leq 1+c_1,\]
\[|h(\lambda)-1|\leq c_1\cdot\lambda^4,\]
where $c_1:= \max_{\lambda\in[0, \frac{1}{2}]}\left(h^{(2)}(\lambda), h^{(4)}(\lambda)      \right).$
\end{lemma}
\noindent\textit{Proof}
The first inequality follows from the estimate
 \[  |h(\tilde{\lambda}) -h(\lambda)| \leq \max _{\lambda_0\in[\lambda, \tilde{\lambda}]} |h^{'}(\lambda_0)|  |\tilde{\lambda}- \lambda|  =  \max _{\lambda_0\in[\lambda, \tilde{\lambda}]} | h^{'}(\lambda_0)-   h'(0)| |\tilde{\lambda}- \lambda|\] 
 
\[\leq\sup_{\tilde{\lambda}_0\in[0, \lambda_0]} |h{''}(\tilde{\lambda}_0)| \cdot \lambda_0 \cdot|\tilde{\lambda}-\lambda| \leq\sup_{\tilde{\lambda}_0\in[0, \frac{1}{2}]} |h{''}(\tilde{\lambda}_0)| \cdot | \tilde{\lambda}+\lambda| \cdot |\tilde{\lambda}-\lambda| \leq c_1 |\tilde{\lambda}^2-\lambda^2    |. \]
 Substituting $\lambda=0$ we get the second inequality. Finally, the third upper bound is obtained using the mean value theorem 4 times.\hfill\(\Box\)

\begin{figure}[h] 
 	\begin{center}
 		\begin{minipage}{5.4cm} 
 			\includegraphics[width=\textwidth]{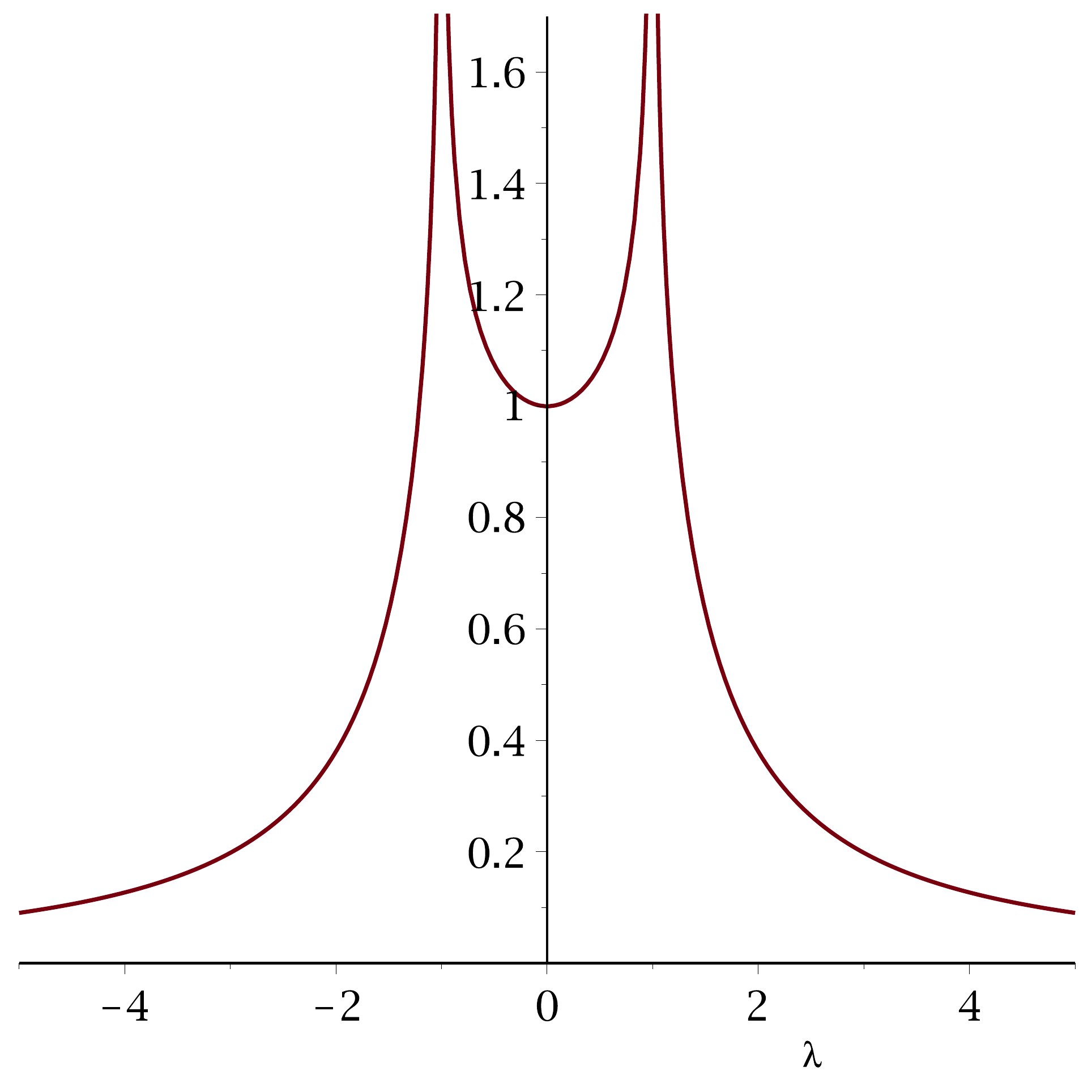}
 		\end{minipage}%
 		\quad \quad
 		\begin{minipage}{5.4cm}
 			\includegraphics[width=\textwidth]{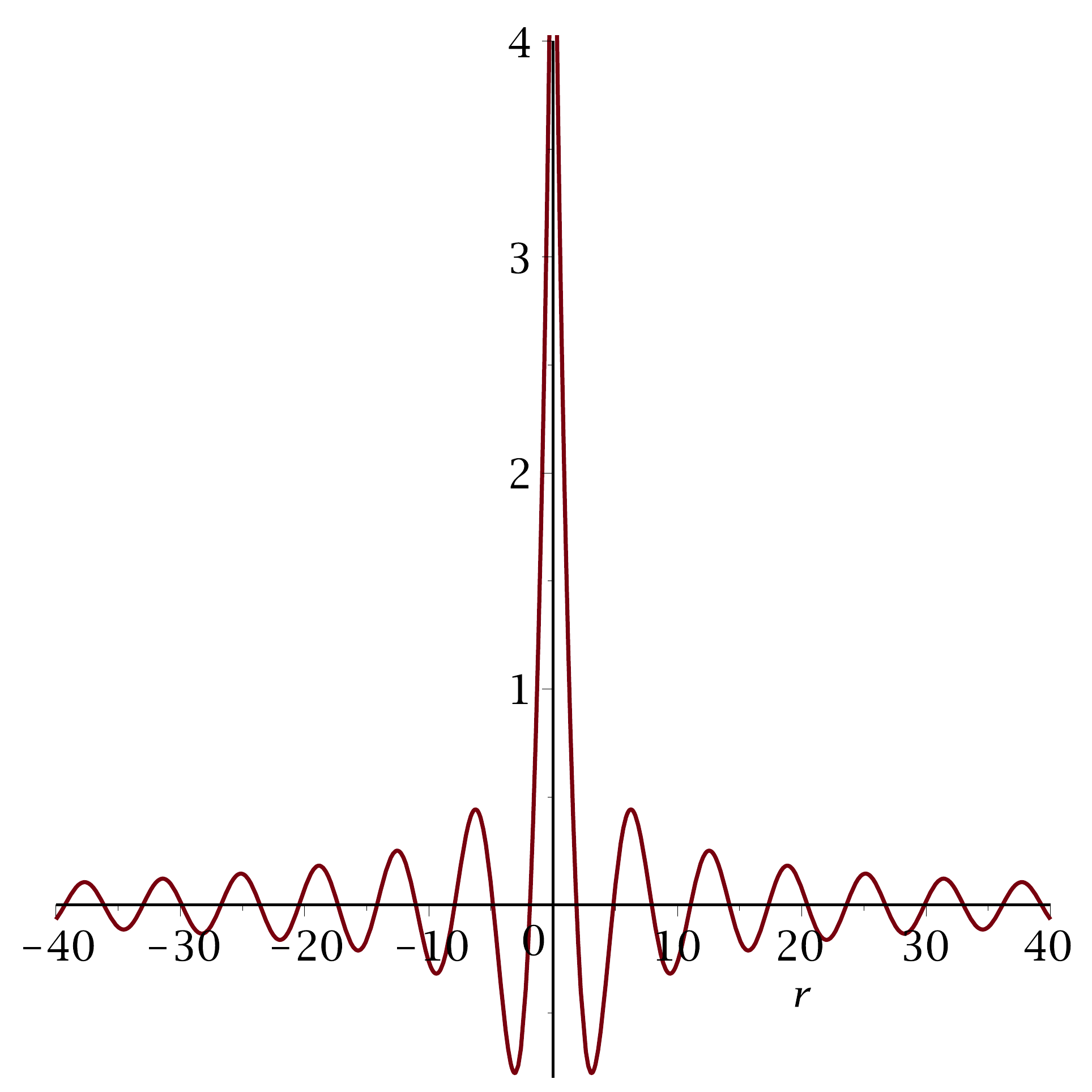}
 		\end{minipage}
 		       \caption{Plots of $f(\lambda)$ and $B(r).$ } \label{fig1}
 	\end{center}
 \end{figure}

\begin{example}\label{remark_2}
The asymptotic behavior of the function $h(\lambda)$ 
when $\lambda \to  \infty$ must guarantee that $f(\cdot)\in L_{1}(\mathbb{R}).$

For example, the function $f\left(\lambda\right)=\left|\lambda^2- s_{0}^2\right|^{-2\alpha} I_{[-M,M]} \left(\lambda\right),$ satisfies Assumption $ \ref{Assumption_1},$
where $~M>s_{0},$ $ I_{[-M,M]}
\left(\lambda\right)= 
\begin{cases} 
1, &\lambda\in [-M,M]\\
0, &\lambda\notin [-M,M]  
\end{cases}$
is the indicator function of the interval $[-M,M].$ 

Another example satisfying Assumption $ \ref{Assumption_1}$ is the following spectral density and corresponding covariance function \[f\left(\lambda\right)= \begin{cases} \frac{1}{\left(1-\lambda^2\right)^{1/4}}, &  \left|\lambda\right|\leq1\\
\frac{1} {\lambda\left(\lambda^2-1\right)^{1/4}}, &  \left|\lambda\right|>1,
\end{cases}
\]
\[B(r)=\sqrt{\pi }\left(\sqrt{2\pi}+\left(\frac{2}{r}\right)^{1/4}\Gamma\left( \frac{3}{4} \right) J_{\frac{1}{4}}\left(r\right)-2  \sqrt {2r}\ {{}_1 F_2}\left(\frac{1}{4};\,\frac{3}{4},\frac{5}{4};\,-\frac{{r}^{2}}{4}\right)\right)\]
 where $J_{\nu}$ is the Bessel function of the first kind, ${{}_1 F_2}$  is the hypergeometric function, and $s_{0}=1$ and $ \alpha=1/8$ were chosen. Plots of  $f\left(\lambda\right)$ and $B(r)$ are shown in Figure~\rm\ref{fig1}. 
 \end{example}
 
\begin{remark} As we study seasonal or cyclic long memory models, in this paper we consider the case of singularities at non-zero frequencies.  The discussion about differences between the cases with spectral singularities at the origin and at other locations can be found in {\rm\cite{ArtRob:1999}.}  As $s_0$ is separated from zero, without loss of generality we assume that $s_{0}>1.$  Indeed, if a time series has a periodic component  with the period $T$ then the corresponding frequency $s_{0}=1/T.$ Changing the time unit the parameter $s_{0}$ can be made greater than 1.
\end{remark}
 
 Now we introduce filter transforms of stochastic processes. To define filters we use real-valued functions $\psi (t),t \in \mathbb{R},$
with the Fourier transforms $\widehat \psi(\cdot).$

Throughout the article, we use the convention that the Fourier transform of an arbitrary function $\psi$ belonging to $L^1(\mathbb{R})$ is the function $\hat{\psi}$ defined, for every $\lambda \in \mathbb{R},$ as \[\hat{\psi}(\lambda)=\int_\mathbb{R}e^{-i\lambda t} \psi(t) dt.\]

\begin{assumption}\label{Assumption_2}
Let $\psi\in L_2 ({\mathbb R})$ be a real-valued function such that $ {\rm supp} \, \widehat \psi\subset [-A,A], A>0,$ and
 $\widehat \psi(\cdot)$ is  continuous except at a finite number of points and of bounded variation on $[-A,A].$
 \end{assumption}
\begin{remark}\label{remark_3}
It follows from Assumption $\ref{Assumption_2}$ 
that  $\widehat \psi$ is bounded and $\psi\in L_{\infty}(\mathbb{R})$ is
 an analytic function.
\end{remark}

Let us define the following constants
$ c_{2}:=  \int_{\mathbb{R}} |\widehat\psi(\lambda)|^{2} d\lambda$ and $  c_{3}:= 2 \int_{\mathbb{R}}\lambda^{2}|\widehat\psi(\lambda)|^{2} d\lambda.$

 Some important for applications functions $\psi (\cdot)$ satisfying 
Assumption \ref{Assumption_2} are the wavelets, see  {\rm\cite{Dau92}, \cite{Meyer92}}, given in the next
 examples. However, in general, $\psi (\cdot)$ is not required to be a wavelet.
\begin{example}\label{example_1}
The function $\psi (\cdot)$ can be selected as the Shannon father or 
mother wavelets. Indeed, the Shannon father wavelet  
\[ \psi_{f} (t)= sinc(\pi t) =
\begin{cases} 
\frac{sin(\pi t)}{\pi t},     &t\neq0,
\\
1,    &t=0,         
 \end{cases}
\] 
has the Fourier transform \[\widehat \psi_{f} (\lambda) = I_{[-\pi,\pi]}(\lambda).\]
 The corresponding  constants  are $ c_{2}=2\pi$  and  $ c_{3}=\frac{4}{3}\pi^3.$

The Shannon mother wavelet \[ \psi_{m} (t)=\frac{sin(2\pi t)- cos(\pi t)}{\pi(t-\frac{1}{2})}
  \]
has the  Fourier transform \[ \widehat \psi_{m} (\lambda) =-\mathrm{e}^{-\frac{i\lambda}{2}} I_{[-2\pi,-\pi)\cup(\pi,2\pi]}(\lambda).\]
  The corresponding  constants are $ c_{2}=2\pi$  and  $ c_{3}=9\frac{1}{3}\pi^3.$

Plots of $\psi_{f} (t)$ and $\psi_{m} (t)$ are shown in  Figure~\rm{\ref{fig22}}.
 \begin{figure}[h]
	\begin{center}
		\begin{minipage}{5.4cm} 
			\includegraphics[width=\textwidth]{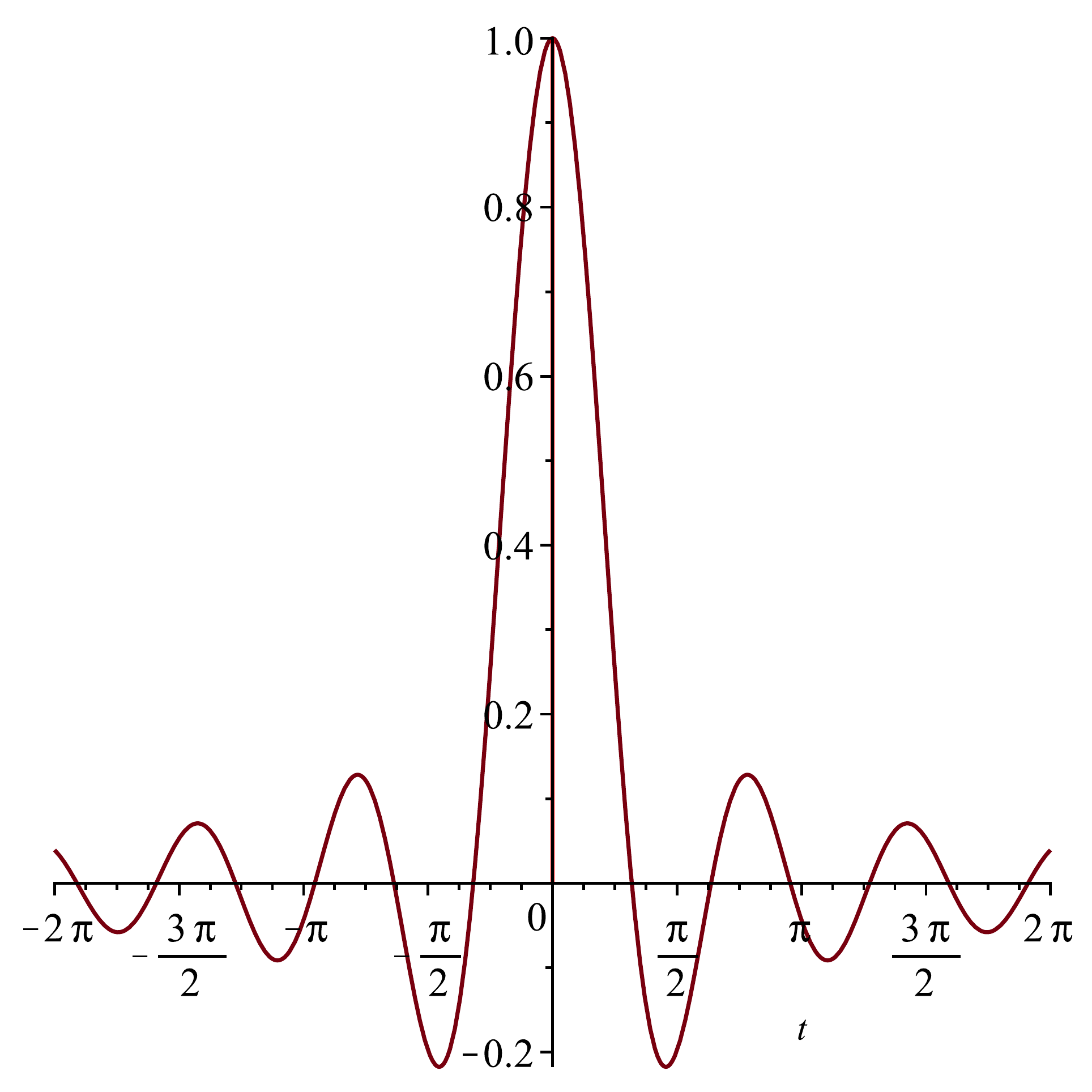}
		\end{minipage}%
		\quad \quad
		\begin{minipage}{5.4cm} 
			\includegraphics[width=\textwidth]{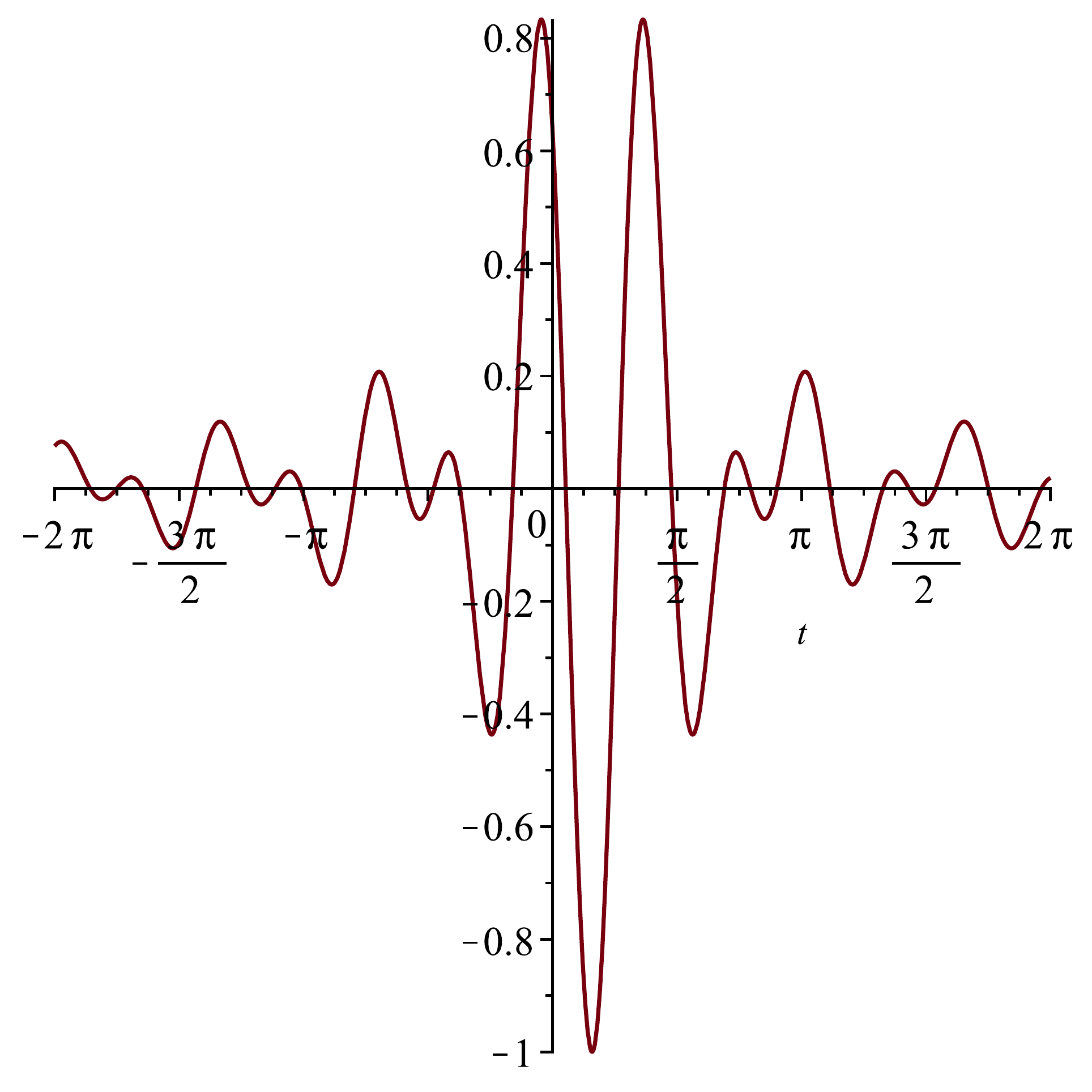}
		\end{minipage}%
		\caption{Plots of $ \psi_{f} (t)$ and $\psi_{m} (t).$} \label{fig22}
	\end{center}
\end{figure}
 \end{example}
\begin{example}\label{example_2}
The function $\psi(\cdot)$ can be selected as the Meyer father 
or mother wavelets (see, for example, {\rm \cite{Dau92},~\cite{Meyer92}}). Indeed, the Meyer wavelets have the Fourier transforms
\[\widehat\psi_{m} (\lambda):= 
\begin{cases} 
  sin\Big(\frac{\pi}{2}\nu\Big(\frac{3|\lambda|}{2\pi}-1\Big)\Big)\mathrm{e}^{\frac{i\lambda}{2}},&   \quad \text{if } \frac{2\pi}{3}\leq|\lambda|\leq\frac{4\pi}{3},
\\
 cos \Big(\frac{\pi}{2}\nu\Big(\frac{3|\lambda|}{4\pi}-1\Big)\Big)\mathrm{e}^{\frac{i\lambda}{2}},&   \quad \text{if } \frac{4\pi}{3}\leq|\lambda|\leq\frac{8\pi}{3},
\\
0, &   \quad \text{otherwise, } 
\end{cases}
\]
\[\widehat\psi_{f} (\omega):= 
\begin{cases} 
1, & \quad \text{if }|\lambda|\leq\frac{2\pi}{3},
\\
  cos \Big(\frac{\pi}{2}\nu\Big(\frac{3|\lambda|}{4\pi}-1\Big)\Big), &   \quad \text{if } \frac{2\pi}{3}\leq|\lambda|\leq\frac{4\pi}{3},
\\
0, &   \quad \text{otherwise, }
 \end{cases}
\]
where $\nu(\cdot)$ is a function with values in [0,1] that satisfies $\nu(x)+\nu(1-x)=1, x\in \mathbb{R}.$ For example, one can use \[ \nu(x)=
\begin{cases}
0, & \quad \text{if } x<0,
\\
x, & \quad \text{if } 0\leq x \leq1,
\\
1, & \quad \text{if } x>1.
 \end{cases}
\]
The corresponding constants for $\widehat\psi_{m} (\lambda)$ are $c_2 =  2\pi  $ and $c_3=12\frac{4}{9}\pi(\pi^2+2)$ and   $c_2 = 2\pi $ and $c_3=1\frac{7}{9}\pi(\pi^2-2)$ for $\widehat\psi_{f} (\lambda)$.

 Plots of $|\widehat\psi_{m} (\lambda)|$ and $ \widehat\psi_{f} (\lambda)$ are shown in {Figure~\rm\ref{fig4}}.
\begin{figure}[h] 
 	\begin{center}
 		\begin{minipage}{5.4cm} 
 			\includegraphics[width=\textwidth]{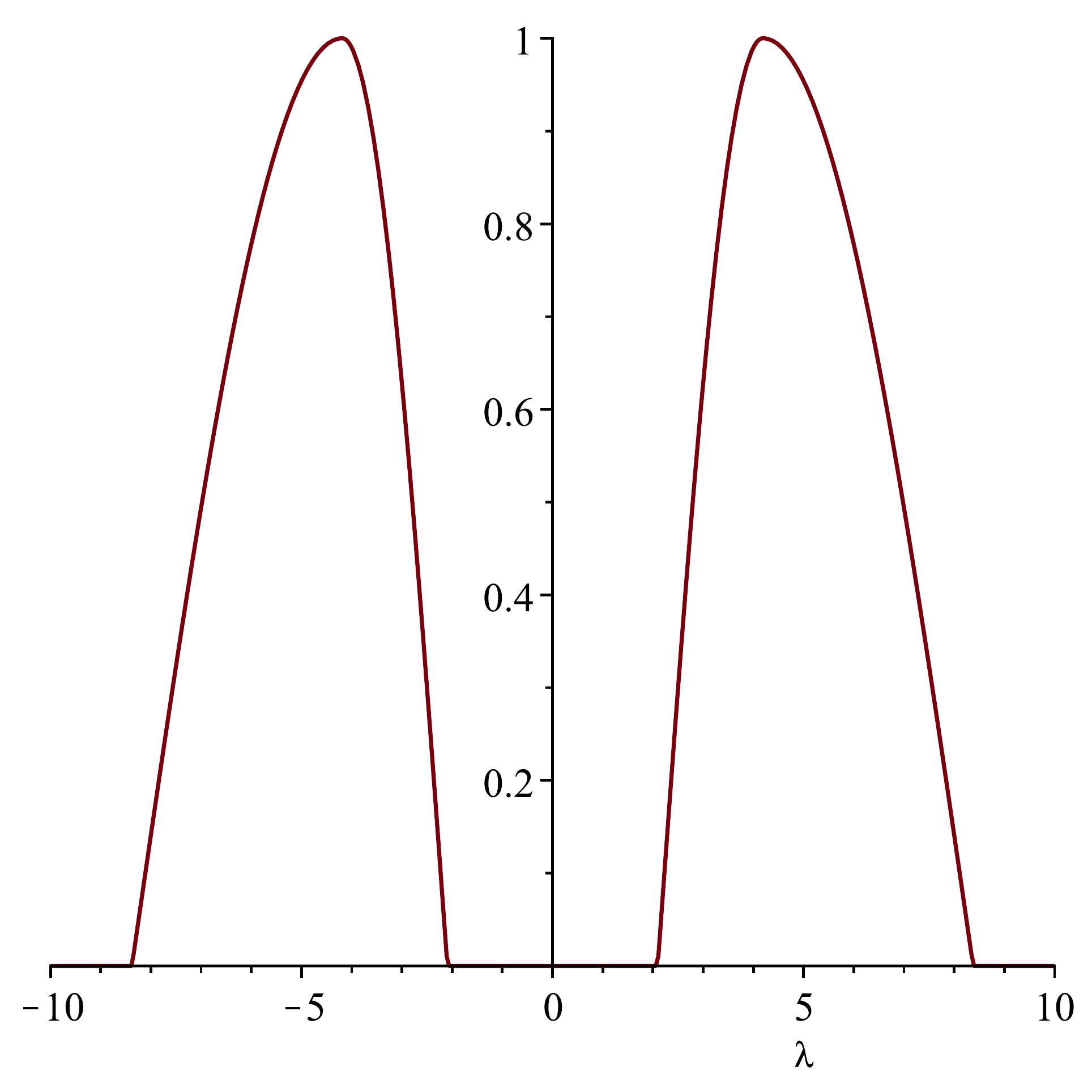}
 		\end{minipage}%
 		\quad \quad
 		\begin{minipage}{5.4cm}
 			\includegraphics[width=\textwidth]{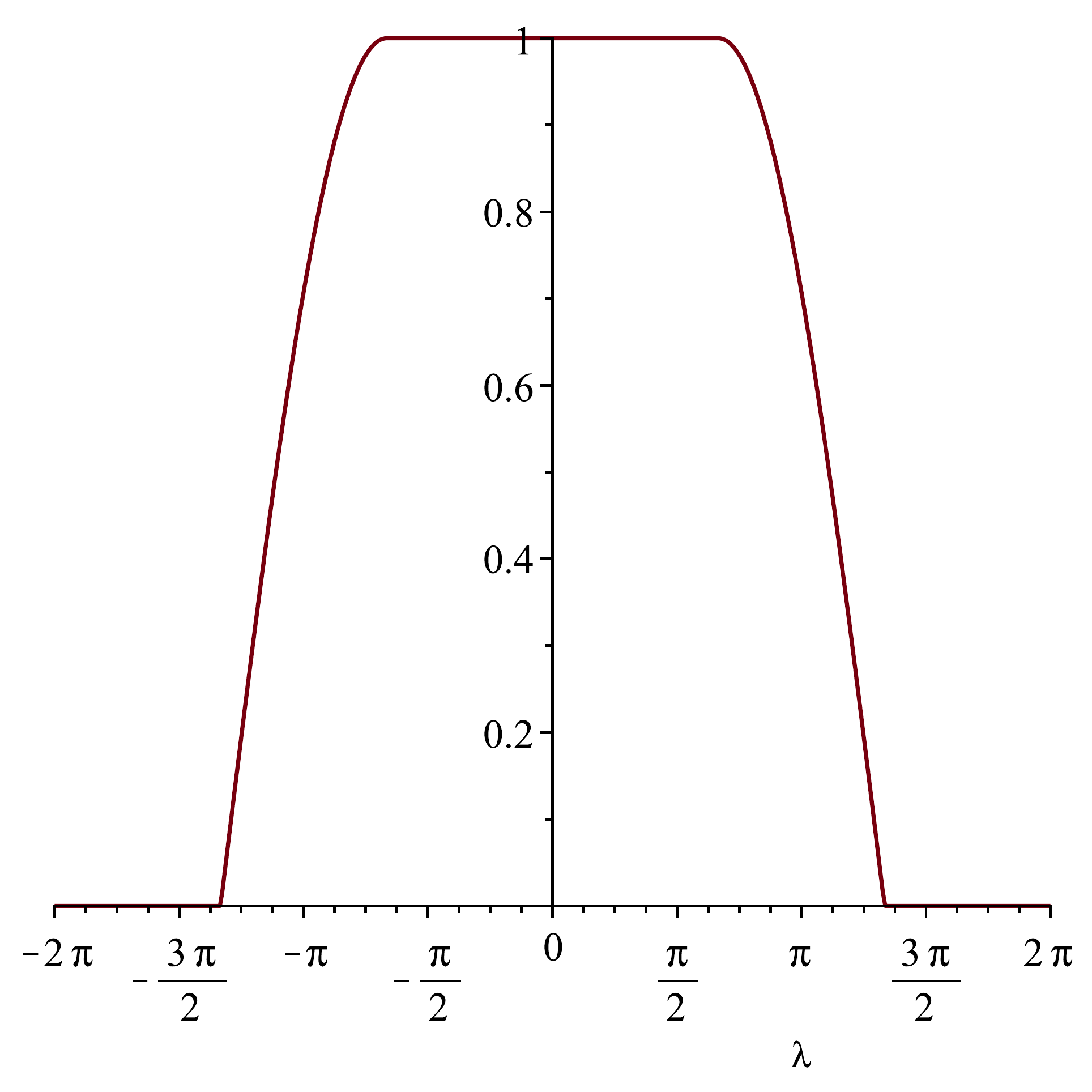}
 		\end{minipage}
 		       \caption{Plots of $|\widehat\psi_{m} (\lambda)|$ and $\widehat\psi_{f} (\lambda).$ } \label{fig4}
 	\end{center}
 \end{figure}
\end{example}
Now, for any pair $(a,b)\in \mathbb{R}_{+}\times \mathbb{R},$ where $\mathbb{R}_{+}=(0,+\infty),$ we define the following filter transform of the process $X(t)$
\begin{align}\label{2_1}
d_{x}(a,b):=\frac{1}{\sqrt{a}} \int_{\mathbb{R}} \psi\bigg(\frac{t-b}{a}\bigg)X(t) dt.
\end{align}
\begin{remark}\label{remark_4}
If $\psi(\cdot)$ is a wavelet, then $d_{x}(a,b) $ given by $(\ref{2_1})$ defines the wavelet transform of the process $X(t)$.\end{remark}
The general filtration theory of stochastic processes guarantees that $(\ref{2_1})$ is correctly defined if the following assumption is satisfied, see  Chapter \rom{5}, \S 6  in \cite{GikSko:2004}.
\begin{assumption}\label{Assumption_3}
Let the integral \[ \int_{\mathbb{R}^2} \psi(t)B(t-t')\overline{\psi(t')}dt dt'\]
exists as an improper Cauchy integral on the plane.
\begin{remark}\label{remark_4'}
Different assumptions on the process and the filter were used in \rm{\cite{Bardet:2010}}. \it{Namely, they assumed that $\psi(\cdot)$ is a mother wavelet that has two vanishing moments and there are constants $c_\psi, c'_\psi>0$ such that $(1+|t|)|\psi(t)| \leq c_{\psi},  |\widehat{\psi}(\lambda)|+|\widehat{\psi}^{\prime}(\lambda)|\leq c'_\psi,$ for all $t,\lambda \in \mathbb{R}.$}
\end{remark}
\end{assumption}
Using the above notations $d_{x}(a,b)$ can be rewritten in the frequency domain as 
\[ d_{x}(a,b)=\sqrt{a}\int_{\mathbb{R}} \mathrm{e}^{i b\lambda}\overline{\widehat{\psi}(a\lambda)}\sqrt{f(\lambda)} dW(\lambda).\]
This Gaussian random variable has a zero mean, i.e. $Ed_{x}(a,b)=0.$ Its variance equals 
\begin{align}\label{2_2}
E|d_{x}(a,b)|^{2}=a  \int_{\mathbb{R}}|\widehat{\psi}(a\lambda)|^{2}f(\lambda) d\lambda:=J(a)
\end{align}
and thus does not depend on $b.$

In the following sections we assume that Assumptions \ref{Assumption_1}-\ref{Assumption_3} are satisfied.

\section{First statistics }\label{sec_3}

Spectral densities satisfying Assumption \ref{Assumption_1} have two parameters of interest $(\alpha$ and $s_{0}).$ This section derives some properties of $d_{x}(a,b)$ and suggests a statistic based on $d_{x}(a,b)$ that can be used as an estimate of $s_{0}^{-4 \alpha}.$ 

Let $\{a_{j}\} ,\{\gamma_{j}\} \subset \mathbb{R_{+}}, \{m_{j}\} \subset \mathbb{N}, j \in \mathbb{N},$ be sequences of positive numbers and $\{b _{jk}\}\subset \mathbb{R}, \, j \in \mathbb{N}, k \in \mathbb{Z},$ be an infinite array. In the following proofs we assume that $\{a_{j}\}$ is an unboundedly monotone increasing sequence and $b_{jk_1}\neq b_{jk_2}$ for all $j\in \mathbb{N}$ and $k_1\neq k_2.$

We will use the following notation
\[
\delta_{jk}:= d_{x}(a_{j},b_{jk})=\sqrt{a_{j}}\int_{\mathbb{R}} \mathrm{e}^{i b_{jk}\lambda}\overline{\widehat{\psi}(a_{j}\lambda)}\sqrt{f(\lambda)} dW(\lambda).
\]
By Assumption \ref{Assumption_1}, 
\[\delta_{jk}= \sqrt{a_{j}}\int_{\mathbb{R}} \mathrm{e}^{i b_{jk}\lambda} \frac{\overline{\widehat{\psi}(a_{j}\lambda)} \sqrt{h(\lambda)}}  {|\lambda^{2}-s_{0}^{2}|^{\alpha}} dW(\lambda).\]
Therefore, for all $j \in\mathbb{N}$ and $ k_1, k_2 \in\mathbb{Z}$ it holds $E \delta_{jk_1}=0 $ and 
\begin{align} \label{3_1}
 I(j,k_{1}, k_{2})&:= Cov(\delta_{j k_{1} },\delta_{j k_{2} })= E(\delta_{j k_{1}}\overline{\delta}_{jk_{2}})\\
                        & =a_{j}\int_{\mathbb{R}} \mathrm{e}^{i (b_{jk_{1} }-b_{jk_{2} }) \lambda}  \frac{|{\widehat{\psi}(a_{j}\lambda)}|^{2} h(\lambda)}  {|\lambda^{2}-s_{0}^{2}|^{2\alpha}} d\lambda\\   
                        &=s_{0}^{-4\alpha}\int_{\mathbb{R}} \mathrm{e}^{i \frac{ b_{jk_{1} }-b_{j k_{2} } }{a_{j}}\lambda}  
\frac {|{\widehat{\psi}(\lambda)}|^{2} h\big(\frac{\lambda}{a_{j}}\big)} {|(\frac{\lambda}{a_{j }s_{0}})^{2}-1|^{2\alpha}} d\lambda.  
\end{align}
\begin{lemma} \label{lemma_1}
Let Assumption \rm\ref{Assumption_2} \it{hold true. Then, for all $k_{1},k_{2}\in \mathbb{Z}^{2}$ and $j \in \mathbb{N}$ such that $\frac{a_{j}}{2A}\geq1$
\[
|I(j, k_{1}, k_{2})|	\leq       c_{4}(s_{0},\alpha)
 \begin{cases}
 1, & \quad \mbox{if } k_{1} =k_{2},\\
 \frac{a_{j}}{|b_{jk_{1}}-b_{jk_{2}}|}, & \quad \mbox{if } k_{1} \neq k_{2},
  \end{cases}
\]
where \[c_{4}( s_{0},\alpha) := 2s_{0}^{-4\alpha} \max\left(\frac{2}{3}c_{2}\left(1+c_{1}\right),   \max_{\lambda\in[-A,A]} \left|\widehat{\psi}(\lambda)\right|^{2} V_{-1/2} ^{1/2} \left(h\left(\cdot\right)\right) 
+\left(1+c_{1}\right) \cdot V_{-A} ^{A} \left( \left|\widehat{\psi}(\cdot)\right|^{2}\right) \right), \]
$V_{a}^{b}(f)$ denotes the total variation of a function $f(\cdot)$ on an interval $[a,b].$}
\end{lemma}
\noindent\textit{Proof} 
If $k_{1}=k_{2}$ then by (\ref{3_1}) 
\[I(j, k_{1}, k_{1})= s_{0}^{-4\alpha}\mathlarger {\int}\limits_{-A}^{A}\frac {|{\widehat{\psi}(\lambda)}|^{2} h\big(\frac{\lambda}{a_{j}}\big)} {\Big|\Big(\frac{\lambda}{a_{j }s_{0}}\Big)^{2}-1\Big|^{2\alpha}} d\lambda \leq  c_{2} s_{0}^{-4\alpha}\sup_{u\in[0,\frac{A}{a_{j}}]} h(u)\bigg|1-\bigg(\frac{A}{a_{j}s_{0}}\bigg)^{2}\bigg|^{-2\alpha}.
 \] 
By the conditions $\frac{a_{j}}{2A}\geq 1$ and $s_0>1$ we get
\begin{align}  \label{3-2}\frac{A}{a_{j}}\leq \frac{1}{2}\quad \textrm{and} \quad \frac{A}{a_{j}s_{0}}\leq \frac{1}{2}.\end{align}
Hence, by Lemma \ref{remark_0} 
\begin{align} \label{3-3}
|I(j, k_{1}, k_{1})|\leq c_{2} s_{0}^{-4\alpha} \frac{1+c_{1}}{(1-\frac{1}{4})^{2\alpha}} \leq\frac{4}{3}(1+c_{1}) c_{2} s_{0}^{-4\alpha}.
\end{align}
 
Let us denote $p\left(\lambda\right):= \frac{\left|\widehat{\psi}\left(\lambda\right)\right|^{2} h\left(\frac{\lambda}{ a_{j}}   \right)} {\left|\left(\frac{\lambda}{a_j s_{0}}\right)^{2} -1 \right|^{2\alpha}}.$ By Assumptions \ref{Assumption_1} and \ref{Assumption_2} the function $p(\cdot)$  is a non-negative integrable function on $[-A,A].$ Therefore, $\tilde{p}\left(\lambda\right):=p\left(\lambda\right)/  \int_{-A}^{A} p(\lambda)d\lambda$
is a probability density. Moreover, by  Assumption \ref{Assumption_1}
we get $ V_{-\frac{1}{2}} ^{\frac{1}{2}} \left( h(\cdot)    \right)=\int_{-\frac{1}{2}}^ {\frac{1}{2}} |h^{'}(\lambda)|d\lambda \leq2 \int_{0}^ {\frac{1}{2}} c_1\cdot \lambda d\lambda<+\infty.$ Hence, it follows from Assumption \ref{Assumption_2}  that 
$\tilde{p}\left(\lambda\right)$ is a function of bounded variation on $[-A,A].$

Therefore, for $k_1\neq k_2$ by Theorem 2.5.3 in \cite{Ushakov:1999} \[\left|\int_{-A}^{A}  \mathrm{e}^{i \frac{ b_{jk_{1} }-b_{j k_{2} } }{a_{j}}\lambda} \tilde{p} \left(\lambda\right) d\lambda \right| \leq \frac{a_j} {\left|b_{jk_{1} }-b_{j k_{2} }  \right| } V_{-A} ^{A} \left( \tilde{p}\right).\]

Hence, if $k_1\neq k_2$ we obtain \[ \left| I\left(j, k_1, k_2\right)\right| \leq s_{0}^{-4\alpha} \int_{-A}^{A} p\left(\lambda\right) d\lambda \cdot V_{-A} ^{A} \left( \tilde{p}\right) \cdot  \frac{a_j} {\left|b_{jk_{1} }-b_{j k_{2} }  \right| } = s_{0}^{-4\alpha} V_{-A} ^{A} \left(p\right) \cdot      \frac{a_j} {\left|b_{jk_{1} }-b_{j k_{2} }  \right| }. \]

It follows from (\ref{3-2}) and $\alpha \in (0, 1/2)$ that
\begin{align*}
V_{-A}^{A} \left(p\right) &\leq 2V_{-A} ^{A} \left(\left| \widehat{\psi}(\cdot)\right|^{2}  h\left(\frac{\cdot}{a_j}\right) \right)\\
&\leq 2\left[\max_{\lambda\in[-A,A]} \left| \widehat{\psi}(\lambda)\right|^{2}  V_{-A} ^{A} \left( h\left ( \frac{\cdot}{a_j} \right) \right) + \max_{\lambda\in[-A,A]}   h\left(\frac{\lambda}{a_j}\right)  
V_{-A} ^{A}  \left( \left| \widehat{\psi}(\cdot)\right|^{2} \right)  \right] \\
&\leq 2\left[ \max_{\lambda\in[-A,A]} \left| \widehat{\psi}(\lambda)\right|^{2} V_{-1/2} ^{1/2}\left( h(\cdot)\right) + \max_{\lambda\in[-\frac{1}{2}, \frac{1}{2}]}h\left(\lambda \right)   V_{-A} ^{A}  \left( \left| \widehat{\psi}(\cdot)\right|^{2} \right) \right].
\end{align*}  
Note that this upper bound does not depend on $ j, a_j, b_{jk}.$

Hence,  \begin{align}\label{3-4} \left| I\left(j, k_1, k_2\right)\right| \leq 2s_{0}^{-4\alpha}  \left[\max_{\lambda\in[-A,A]} \left| \widehat{\psi}(\lambda)\right|^{2}    V_{-1/2} ^{1/2} \left(h\left(\cdot \right) \right)   +  \max_{\lambda\in[-\frac{1}{2}, \frac{1}{2}]}   h\left(\lambda \right) \cdot   V_{-A} ^{A} \left( \left| \widehat{\psi}(\cdot)\right|^{2} \right) \right] \end{align}

Comparing (\ref{3-3}) and (\ref{3-4}) we obtain the value of $ c_{4}(s_{0},\alpha)$ at the statement of Lemma \ref{lemma_1}.\hfill\(\Box\)
  
 Let us define
  \[\bar{\delta}_{j\cdot }^{(2)}:=\frac{1}{m_{j}} \displaystyle\sum_{k=1}^{m_{j}} \delta_{jk}^{2}.\]
 It follows from (\ref{2_2}) that $ E\bar{\delta}_{ j\cdot}^{(2)}= J(a_{j}).$   
\begin{lemma} \label{lemma_2}
Suppose that $\frac{a_{j}}{2A}\geq 1,$ the sequences $\{b_{jk}\}$ and $\{\gamma_{j}\}$ are such that  for all $ j\in \mathbb{N}, k_{1}, k_{2}\in\mathbb{Z},$ it holds 
$|b_{jk_{1}}-b_{j k_{2}}|\geq|k_{1}-k_{2}| \gamma_{j}.$ Then,
 \[ 
Var\Big(\bar{\delta}_{ j\cdot}^{(2)}\Big)\leq \frac{ c_{5j} ( s_{0},\alpha) }{m_{j}},
\]
where $ c_{5j} ( s_{0},\alpha):= 2 c_{4}^{2}\Big( s_{0},\alpha \Big) \Big(1+\frac{\pi^{2}}{3}\frac{a_{j}^{2}}{\gamma_{j}^{2}}\Big).$
 \end{lemma} 
\noindent\textit{Proof}
 Notice, that by (\ref{2_1})  random variables $\delta _{jk}$ are  centered Gaussian. For any centered 2-di\-mensional Gaussian vector $(z_{1}, z_{2})$ it holds $Cov(z_{1}^{2}, z_{2}^{2}) = 2 Cov^{2}(z_{1}, z_{2}).$ Therefore, recalling (\ref{2_2}) and the definition of $I(j, k_{1}, k_{2})$ in (\ref{3_1}), we obtain
\begin{align*}
Var\Big(\bar{\delta}_{ j \cdot}^{(2)}\Big)&= Var\bigg(\frac{1}{m_{j}} \displaystyle\sum_{k=1}^{m_{j}} \delta_{jk}^{2}\bigg)=\frac{1}{m_{j}^{2}}\displaystyle\sum_{1\leq k_{1}, k_{2}\leq m_{j}}Cov \Big(\delta_{ jk_{1} }^{2}, \delta_{ jk_{2} }^{2}\Big)\\
&=\frac{2}{m_{j}^{2}}\displaystyle\sum_{1\leq k_{1}, k_{2}\leq m_{j}}I^{2}\left (j, k_{1}, k_{2}\right)=\frac{2}{m_{j}^{2}}\left(\displaystyle\sum_{ k_{1}=1}^{m_{j}}I^{2}\left(j, k_{1}, k_{1}\right)
 + \sum_{\substack{1\leq k_{1}, k_{2}\leq m_{j}\\ k_{1}\neq k_{2}}} I^{2}\left(j, k_{1}, k_{2}\right) \right).
\end{align*}
By Lemma \ref{lemma_1} it follows that 
\begin{align*}
Var\left(\bar{\delta}_{ j\cdot}^{(2)}\right)&\leq \frac{2}{m_{j}^{2}} \left(c_{4}^{2}(s_{0},\alpha )m_{j} +  2 c_{4}^{2}(s_{0},\alpha )\displaystyle\sum_{k_{1}=1}^{m_{j}}  \displaystyle\sum_{k_{2}=k_{1}+1}^{m_{j}} \frac {a^{2}_{j}} {\left(b_{jk_{1}}- b_{jk_{2}}\right)^2}\right)\\&\leq
\frac{2      c_{4}^{2}(s_{0},\alpha )} {m^{2}_{j}} \left( m_{j}+\frac{2a^{2}_{j}} {\gamma^{2}_{j}}
\displaystyle\sum_{ k_{1}=1}^{m_{j}}\ \displaystyle\sum_{ k_{2}=k_{1}+1}^{m_{j}} \frac{1}{\left(k_{2}-k_{1}\right)^{2}}\right)\\ &\leq \frac{2       c_{4}^{2}(s_{0},\alpha )}{m^{2}_{j}} \left(m_{j} + \frac{2a^{2}_{j}} {\gamma^{2}_{j}} \displaystyle\sum_{ k_{1}=1}^{m_{j}}\ \displaystyle\sum_{ \tilde{k}_{2}=1}^{\infty}  \frac{1}{\tilde{k}_{2}^{2}} \right)=\frac{2       c_{4}^{2}(s_{0},\alpha )} {m_{j}} \left(1+\frac{\pi^{2}}{3}\frac{a_{j}^{2}}{\gamma_{j}^{2}}\right).      
\end{align*}\hfill\(\Box\)
 \begin{lemma} \label{lemma_3}
Let $\{r_{j}\}\subset \mathbb{R_{+}}$ be a decreasing sequence such that $\displaystyle\lim_{j \to \infty} r_{j}=0.$ Let us choose such $\{m_{j}\}$ that  $\displaystyle\sum_{j=1}^{\infty} \frac{1}{r^{2}_{j}m_{j}}<+\infty$ and $\displaystyle\sum_{j=1}^{\infty} \frac{a_{j}^{2}} {r_{j}^{2} \gamma_{j}^{2} m_{j} }<+\infty, $ where $\{ \gamma_{j}\}$ is from Lemma \rm {\ref{lemma_2}}. \it{Then there exists an almost surely finite random variable $ c_{6}$ such that for all $j\in \mathbb{N}$
\[\left|\bar{\delta}_{ j\cdot}^{(2)} -J(a_{j})\right|\leq c_{6}r_{j}.\]}
 \end{lemma}

\noindent\textit{Proof}
Using Chebyshev inequality and Lemma \rm {\ref{lemma_2}} we obtain 
\[P\left(\left|\bar{\delta}_{ j\cdot}^{(2)} -J(a_{j})\right|>r_{_{j}}\right) \leq \frac{Var(\bar{\delta}_{ j\cdot}^{(2)})}{r_{j}^{2}}\leq \frac{  c_{5j}(s_{0},\alpha )}{ r_{j}^{2}m_{j}}= 2      c_{4}^{2}(s_{0},\alpha )\frac{\Big(1+\frac{\pi^{2}}{3}\frac{a_{j}^{2}}{\gamma_{j}^{2}}\Big)} { r_{j}^{2}m_{j}}.\]
By the choice of $\{m_{j}\}$ 
 \[\displaystyle\sum_{j=1}^{\infty}P\left(\left|\bar{\delta}_{ j\cdot}^{(2)} -J(a_{j})\right|>r_{_{j}}\right)<+\infty.\] 
Therefore, applying  the Borel–Cantelli lemma we obtain the required statement.\hfill\(\Box\)

We will use the next technical result.
 \begin{lemma} \label{lemma_4}
 For all $ \alpha \in (0, \frac{1}{2}) $ and $ x\in [0, \frac{1}{2}]$ it holds \[0\leq(1-x)^{-2\alpha}-1\leq4x. \]
\end{lemma}
\noindent\textit{Proof}
Applying the mean value theorem to the function $\gamma_{\alpha}(y)=(1-y)^{-2\alpha}$ on the interval $[0,x]$ we obtain
\[ (1-x)^{-2\alpha}-1=\gamma_{\alpha}(x)-\gamma_{\alpha}(0)
=\gamma^{\prime}_{\alpha}(a_{\alpha})x=2\alpha (1- a_{\alpha})^{-2\alpha-1}x \leq   (1- a_{\alpha})^{-2\alpha-1}x, \]  
where $a_{\alpha}\in [0,x]\subset[0, \frac{1}{2}].$  \par

Therefore, as $ \alpha\in (0,\frac{1}{2})$ we get  
$(1-x)^{-2\alpha}-1\leq\Big(1- \frac{1}{2}\Big)^{-2}x=4x. $ \hfill\(\Box\) 

The lemma below gives an upper bound  on the deviation of $ J(a_{j}) $ from $s_{0}^{-4\alpha}c_{2}.$

\begin{lemma} \label{lemma_5}
If $j\in \mathbb{N}$  is such that $\frac{a_{j}}{A}\geq2,$ then one has 
\[\left|J(a_{j})- c_{2} s_{0}^{-4\alpha}\right|\leq\frac{c_{7}(s_{0},\alpha)}{a^2_{j}},\]  
where $  c_{7}(s_{0},\alpha ) := \left(\frac{4(1+c_1)}{s^2_{0}}+c_1\right) c_{2} \frac{A^2}{s^{4\alpha}_{0}}.$
\end{lemma}
\noindent\textit{Proof}
Noting that $J(a_{j})=I (j, k_{1},k_{1})$ and using (\ref{3_1}) we get 
\[J(a_{j})- c_{2} s_{0}^{-4\alpha} =s_{0}^{-4\alpha}  \int_{-A}^{A}\left|\widehat\psi(\lambda)\right|^{2} \left(\frac{h\left(\frac{\lambda}{a_{j}}\right)}{\left|1-\left(\frac{\lambda}{a_{j}s_{0}}\right)^{2}\right|^{2\alpha}}-1\right) d\lambda.\] 
Now, by Lemma \ref{lemma_4}, the conditions on $h(\cdot)$ in Assumption \ref{Assumption_1}  and Lemma \ref{remark_0}  it follows
\[
\bigg|J(a_{j})-  c_{2}s_{0}^{-4\alpha} \bigg|\leq s_{0}^{-4\alpha}  \int_{-A}^{A}\left|\widehat\psi(\lambda)\right|^{2}\Bigg|h\Bigg(\frac{\lambda}{a_{j}}\Bigg)\Bigg(\Bigg(1- \Bigg(\frac{\lambda}{a_{j}s_{0}}\Bigg)^{2}\Bigg)^{-2\alpha} -1 \Bigg)\]
\[ +h\Bigg(\frac{\lambda}{a_{j}}\Bigg)-1 \Bigg|d\lambda  \leq s_{0}^{-4\alpha}  \int_{-A}^{A}\left|\widehat\psi(\lambda)\right|^{2} \left(4h\Bigg(\frac{\lambda}{a_{j}}\Bigg) \Bigg(\frac{\lambda}{a_{j}s_{0}}\Bigg)^{2} + c_1\left(\frac{\lambda}{a_j}\right)^2  \right) d\lambda.
\]  
Moreover, it follows from  Lemma \ref{remark_0} and the conditions of the lemma that for $\lambda\in[-A,A]$
\[
 4h\Bigg(\frac{\lambda}{a_{j}}\Bigg) \Bigg(\frac{\lambda}{a_{j}s_{0}}\Bigg)^{2} +c_1\left(\frac{\lambda}{a_{j}}\right)^2 \leq\left(4 (1+c_1) \frac{A^2}{s^2_{0}} + c_1A^2\right) a^{-2}_j.
\] 
which completes the proof. \hfill\(\Box\) 

Combining Lemmas \ref{lemma_3} and \ref{lemma_5} we obtain
\begin{proposition}\label{pro_1}
Under the conditions of Lemma \rm {\ref{lemma_3}} \it{it holds $\bar{\delta}_{ j\cdot}^{(2)}\xrightarrow{a.s.} c_{2} s_{0}^{-4\alpha}  .$ Moreover, there exists an almost surely finite random variable $ c_{8}$ such that  for all $j\in \mathbb{N}$ 
\[
\Big|\bar{\delta}_{ j\cdot}^{(2)}- c_{2}s_{0}^{-4\alpha} \Big|\leq c_{8}\max(r_{j}, a_{j}^{-2}).
\]}
\end{proposition}
\section{Second statistics}\label{sec_{4}}

In this section we further study properties of $\bar{\delta}_{ j\cdot}^{(2)}$ and $ J(a_{j}).$ It allows us to suggest a new estimate of $ \alpha s_{0}^{-4\alpha -2}.$ The main idea is to find the asymptotic behaviour of increments of $\bar{\delta}_{ j\cdot}^{(2)}.$ Therefore, we start by deriving some results about increments of $J(a_{j})=\mathbb{E}\bar{\delta}_{ j\cdot}^{(2)}.$
\begin{lemma} \label{lemma_6} If $\{a_{j}\}$ is an unboundedly monotone increasing then
 \[\lim_{j\to +\infty} \frac {J(a_{j})-  J(a_{j+1})} {a_{j}^{-2}-a_{j+1}^{-2} } =  \alpha c_{3} s_{0}^{-4\alpha-2}.\]
   \end{lemma}
 \noindent\textit{Proof}
  By (\ref{2_2}) and Assumption \ref{Assumption_1}
we get 
\[
\frac{ J(a_{j+1})-J(a_{j})} {a_{j}^{-2}-a_{j+1}^{-2} }=\int_{\mathbb{R}}\frac{|\widehat\psi(\lambda)|^{2}} {a_{j}^{-2}-a_{j+1}^{-2} } \left(\frac{h\Big(\frac{\lambda}{a_{j+1}}\Big)}{\Big|s_{0}^{2}-\Big(\frac{\lambda}{a_{j+1}}\Big)^{2}\Big|^{2\alpha}} -  \frac{h\Big(\frac{\lambda}{a_{j}}\Big)}{\Big|s_{0}^{2}-\Big(\frac{\lambda}{a_{j}}\Big)^{2}\Big|^{2\alpha}}    \right)  d\lambda\]
 \[= \int_{-A}^{A}|\widehat\psi(\lambda)|^{2} \frac{ h\Big(\frac{\lambda}{a_{j+1}} \Big) \Big|s_{0}^{2}-\Big(\frac{\lambda}{a_{j}}\Big)^{2}\Big|^{2\alpha}- h\Big(\frac{\lambda}{a_{j}} \Big) \Big|s_{0}^{2}-\Big(\frac{\lambda}{a_{j+1}}\Big)^{2}\Big|^{2\alpha} }
{  \left( a_{j}^{-2}-a_{j+1}^{-2}\right)          \Big|s_{0}^{2}-\Big(\frac{\lambda}{a_{j+1}}\Big)^{2} \Big|^{2\alpha} \Big |s_{0}^{2}-\Big(\frac{\lambda}{a_{j}}\Big)^{2}\Big|^{2\alpha}} d\lambda.\]
 As $\lambda\in [-A,A]$ and $a_{j}\to +\infty$ when $ j\to +\infty,$ then there is $ j_{0}\in \mathbb{N},$ such that $\frac{|\lambda|}{a_{j}} \leq \frac{1} {2},$ for all $\lambda\in [-A,A]$ and $j\geq j_{0}.$
Hence, for sufficiently large $j\geq j_{0}$ the integrand can be bounded as it is shown below

$\left|\widehat\psi(\lambda)\right|^{2} \frac{\bigg| h\Big(\frac{\lambda}{a_{j+1}} \Big) \Big|s_{0}^{2}-\Big(\frac{\lambda}{a_{j}}\Big)^{2}\Big|^{2\alpha}- h\Big(\frac{\lambda}{a_{j}} \Big) \Big|s_{0}^{2}- \Big(\frac{\lambda}{a_{j+1}}\Big)^{2}\Big|^{2\alpha} \bigg|}
{     \left(a_{j}^{-2}-a_{j+1}^{-2} \right)         \Big|s_{0}^{2}- \Big(\frac{\lambda}{a_{j+1}}\Big)^{2}\Big|^{2\alpha}\Big|s_{0}^{2}-\Big(\frac{\lambda}{a_{j}}\Big)^{2}\Big|^{2\alpha}}
 $
\begin{align}
 &\leq  \frac{ \left|\widehat\psi(\lambda)\right|^{2}\bigg| h\Big(\frac{\lambda}{a_{j+1}} \Big) \Big|s_{0}^{2}-\Big(\frac{\lambda}{a_{j}}\Big)^{2}\Big|^{2\alpha}- h\Big(\frac{\lambda}{a_{j}} \Big) \Big|s_{0}^{2}- \Big(\frac{\lambda}{a_{j+1}}\Big)^{2}\Big|^{2\alpha} \bigg|} 
{\left(a_{j}^{-2}-a_{j+1}^{-2} \right)     \left|s_0^2-s_0^2\max_{\lambda\in[-A,A] }\left( \frac{\lambda}{a_{j+1}} \right)^2\right|^{2\alpha}   \left|s_0^2-s_0^2\max_{\lambda\in[-A,A] }\left( \frac{\lambda}{a_{j}} \right)^2\right|^{2\alpha}       }\\
&\leq\left|\widehat\psi(\lambda)\right|^{2} \frac{\left| h\left (\frac{\lambda}{a_{j+1}} \right)  - h\left (\frac{\lambda}{a_{j}} \right) \right|  \left| s^2_0 -\left(\frac{\lambda}{a_j}\right)^2\right|^{2\alpha }+h\left(\frac{\lambda}{a_j} \right) \left|   \left|   s_0^2-  \left(\frac{\lambda}{a_j}   \right)^2 \right|^{2\alpha}    -  \left|s_0^2-  \left(\frac{\lambda}{a_{j+1}}   \right)^2 \right|^{2\alpha} \right|  }
 {\left(3s_0^2/4\right)^{4\alpha}  \left(a_{j}^{-2}-a_{j+1}^{-2} \right)}.\label{4_11}  \end{align}
Applying the first inequality in Lemma \ref{remark_0} to $ h\left (\frac{\lambda}{a_{j+1}}\right) -    h\left (\frac{\lambda}{a_{j}}\right)$ and the mean value theorem to $\left(s_0^2- \left (\frac{\lambda}{a_{j}}\right)^2    \right)^{2\alpha}    -    \left(s_0^2- \left (\frac{\lambda}{a_{j+1}}\right)^2    \right)^{2\alpha}$ we obtain that the upper bound for the right part of (\ref{4_11}) is 
\[\left|\widehat\psi(\lambda)\right|^{2} \frac{c_1 s_0^{4\alpha}\lambda^2+2\alpha\left(1+c_1\right)  \left(  3s_0^2/4 \right)^{2\alpha-1}\lambda^2}{   \left(  3s_0^2/4 \right)^{2\alpha}}. \]

This upper bound is integrable and does not depend on $j.$ Hence, one can use the dominated convergence theorem. For $\lambda\in [-A,A]$ it holds $\lim_{j\to+\infty}    h\left (\frac{\lambda}{a_{j}}\right)=1$ and 
$\lim_{j\to+\infty}   \left| s_0^2- \left (\frac{\lambda}{a_{j}}\right)^2    \right|^{2\alpha}=~s_0^{4\alpha}. $
Hence,
\begin{align*}\label{4_1}
  \lim_{j\to +\infty}\frac{ h\Big(\frac{\lambda}{a_{j+1}} \Big) \Big|s_{0}^{2}-\Big(\frac{\lambda}{a_{j}}\Big)^{2}\Big|^{2\alpha}-   h\Big(\frac{\lambda}{a_{j}} \Big) \Big|s_{0}^{2}-\Big(\frac{\lambda}{a_{j+1}}\Big)^{2}\Big|^{2\alpha} }
{  \left(a_{j}^{-2}-a_{j+1}^{-2} \right )    \Big|s_{0}^{2}-\Big(\frac{\lambda}{a_{j+1}}\Big)^{2}\Big|^{2\alpha}   \Big|s_{0}^{2}-\Big(\frac{\lambda}{a_{j}}\Big)^{2}\Big|^{2\alpha}} 
 \end{align*}
 \begin{align}
 &=\lim_{j\to +\infty}\frac{ h\Big(\frac{\lambda}{a_{j+1}}\Big)\Big|s_{0}^{2}- \Big(\frac{\lambda}{a_{j}}\Big)^{2}\Big|^{2\alpha}} {   \left(a_{j}^{-2}-a_{j+1}^{-2} \right ) s_{0}^{8\alpha}} \left(1- \frac{ h\Big(\frac{\lambda}{a_{j}}\Big)} { h\Big(\frac{\lambda}{a_{j+1}}\Big)} \Bigg|  \frac{s_{0}^{2}- \Big(\frac{\lambda}{a_{j+1}}\Big)^2  } {s_{0}^{2}- \Big(\frac{\lambda}{a_{j}}\Big)^2  }\Bigg|^{2\alpha} \right)  \\
&=\lim_{j\to +\infty} \frac{s_{0}^{-4\alpha}} {a_{j}^{-2}-a_{j+1}^{-2} }  \left(1- \frac{ h\Big(\frac{\lambda}{a_{j}}\Big)} { h\Big(\frac{\lambda}{a_{j+1}}\Big)}  \Bigg| 1+ \frac{ \Big(\frac{\lambda}{a_{j}}\Big)^{2} -  \Big(\frac{\lambda}{a_{j+1}}\Big)^{2} } {s_{0}^2 -\Big(\frac{\lambda}{a_{j}}\Big)^{2} }\Bigg|^{2\alpha} \right). 
\end{align} 
Using L'Hopital's rule, one can see that for $\alpha \in (0,1/2)$ it holds $ \lim_{x\to 0}\frac{1-(1+x)^{2\alpha}} {x}=- 2\alpha.$

Noting that for $\lambda \in [-A,A]$ we get  \[ \frac{ \left|h\left(\frac{\lambda}{a_{j+1}}\right)-h\left(\frac{\lambda}{a_{j}}\right)\right  |}          { a_j^{-2}-a_{j+1}^{-2} } \leq \sup_{
\lambda_0\in[0,A]}\left| h^{''} \left( \frac{\lambda_0}{a_j}  \right)  \right| \cdot 
\lambda^2 \to 0,\]   
when $a_j\to +\infty,$ we obtain that (\ref{4_1}) equals
 \begin{align*}
&\lim_{j\to +\infty} \frac{ s_{0}^{-4\alpha}} {\left(a_{j}^{-2}-a_{j+1}^{-2}\right)} \frac{h\Big(\frac{\lambda}{a_{j}}\Big)}
{ h\Big(\frac{\lambda}{a_{j+1}} \Big)} \left( \frac{ h\Big(\frac{\lambda}{a_{j+1}}\Big) } 
{h\Big(\frac{\lambda}{a_{j}}\Big)}  -1+1-\Bigg| 1+  \frac{ \Big(\frac{\lambda}{a_{j}}\Big)^{2} -  \Big(\frac{\lambda}{a_{j+1}}\Big)^{2} } {s_{0}^2 -\Big(\frac{\lambda}{a_{j}}\Big)^{2} }\Bigg|^{2\alpha}\right)\\
&= \lim_{j\to +\infty}  \frac{s_{0}^{-4\alpha}}  {\left(a_{j}^{-2}-a_{j+1}^{-2}\right)} \Bigg( 1-\Bigg| 1+  \frac{ \Big(\frac{\lambda}{a_{j}}\Big)^{2} -  \Big(\frac{\lambda}{a_{j+1}}\Big)^{2} } {s_{0}^2 -\Big(\frac{\lambda}{a_{j}}\Big)^{2} }\Bigg|^{2\alpha}\Bigg)=-\lim_{j\to +\infty} \frac{ 2\alpha s_{0}^{-4\alpha} \left(\Big(\frac{\lambda}{a_{j}}\Big)^{2} -  \Big(\frac{\lambda}{a_{j+1}}\Big)^{2}\right) } { \left(a_{j}^{-2}-a_{j+1}^{-2}\right)       \left(s_{0}^2 -\left(\frac{\lambda}{a_{j}}\right)^{2}\right) }\\
&= - 2\alpha s_{0}^{-4\alpha-2}\lambda^{2}.  \hspace{11.1cm} \Box 
\end{align*}

Now we investigate the rate of convergence in Lemma \ref{lemma_6}.
 \begin{lemma} \label{lemma_7}
 There is $j_0\in \mathbb{N}$ such that for all $j\geq j_0$ it holds
\begin{align*}&\Bigg|\frac {J(a_{j})-  J(a_{j+1})} {a_{j}^{-2}-a_{j+1}^{-2} } -\alpha s_{0}^{-4\alpha-2}  c_{3}\Bigg| \leq \frac{ c_{3}2^{6}A^{2}  s_{0}^{-4\alpha-4} \left(1+ 3^3c_1 s^4_{0}/2^6\right) } {3^3a_{j}^{2}\Big(1- \Big(\frac{a_{j}} {a_{j+1}}\Big)^{2}\Big)}.
\end{align*}
\end{lemma} 
\noindent\textit{Proof}
 Note that 
 \begin{align*}
 &\Bigg|J(a_{j})-  J(a_{j+1})  -\alpha s_{0}^{-4\alpha-2} c_{3}\cdot (a_{j}^{-2}-a_{j+1}^{-2})\Bigg|\\
 &= \Bigg|\int_{-A}^{A}|\widehat\psi(\lambda)|^{2}\Bigg(  \frac{ h\Big(\frac{\lambda}{a_{j}}\Big)} {\Big|s_{0}^2 -\Big(\frac{\lambda}{a_{j}}\Big)^{2}\Big|^{2\alpha} }-2\alpha s_{0}^{-4\alpha-2} \Big(\frac{\lambda}{a_{j}}\Big)^{2}-s_0^{-4\alpha}   \Bigg) d\lambda \\ 
&  -\int_{-A}^{A}|\widehat\psi(\lambda)|^{2}\Bigg(  \frac{ h\Big(\frac{\lambda}{a_{j+1}}\Big)} {\Big|s_{0}^2 -\Big(\frac{\lambda}{a_{j+1}}\Big)^{2}\Big|^{2\alpha} }    -2\alpha s_{0}^{-4\alpha-2} \Big(\frac{\lambda}{a_{j+1}}\Big)^{2} -s_0^{-4\alpha}  \Bigg) d\lambda \Bigg|
\end{align*}
\begin{align}
 &\leq s_{0}^{-4\alpha}\int_{-A}^{A}|\widehat\psi(\lambda)|^{2}\Bigg (\frac{|h \big(\frac{\lambda}{a_{j}}\big)-1|}  {|1- \big(\frac{\lambda}{a_{j}s_{0}}\big)^{2}| ^{2\alpha} }+\bigg| \frac{1}  {{|1- \big(\frac{\lambda}{a_{j}s_{0}}\big)^{2}| ^{2\alpha} }} -2\alpha\bigg( \frac{\lambda}{a_{j}s_{0}}\bigg)^{2} -1\bigg|\\ 
 & +  \frac{|h \big(\frac{\lambda}{a_{j+1}}\big)-1|}  {|1- \big(\frac{\lambda}{a_{j+1}s_{0}}\big)^{2}| ^{2\alpha} }+
\bigg| \frac{1}  {{|1- \big(\frac{\lambda}{a_{j+1}s_{0}}\big)^{2}| ^{2\alpha} }} -2\alpha\bigg( \frac{\lambda}{a_{j+1}s_{0}}\bigg)^{2} -1\bigg|   \Bigg)d\lambda.\label{4_2}
\end{align}
  Let us consider the function $ f(x):=(1-x)^{-2\alpha}-2\alpha x-1, x\in[0,1/4], \alpha \in(0,1/2).$ Notice, that  $ f(\cdot)\in C^{2}[0,1/4], f(0)=0,$ $ f^{\prime}(x)=2\alpha (1-x)^{-2\alpha-1} -2\alpha,$ $f^{\prime}(0)=0,$ $f'' (x) =2 \alpha (2\alpha+1)(1-x)^{-2\alpha-2}.$ Then, applying the mean value theorem twice, we get    
 \[|f(x)-f(0)| \leq \sup_{{u}\in [0,1/4]} f'' (u)x^2, \quad x\in[0,1/4].\] 
   Noting that $\sup_{{u}\in [0,1/4]} f'' (u) \leq \frac{2^7}{3^3}$ it follows 
 $|f(x)-f(0)| \leq  \frac{2^7}{3^3} x^2, \quad x\in[0,1/4]. $
 
 Therefore, if $\frac{A}{a_{j}s_{0}}\leq \frac{1}{2}$ one can bound the second and the fourth terms of the integrand in (\ref{4_2}) by $ \frac{2^7}{3^3} \big(\frac{\lambda}{a_{j}s_{0}}\big)^4 $ and $\frac{2^7}{3^3}\big(\frac{\lambda}{a_{j+1}s_{0}}\big)^4 $ respectively.
  By  Lemma $ \ref{remark_0},$ if $ \frac{A}{a_{j}}\leq\frac{1}{2}$ then the first and third terms in (\ref{4_2}) can be bounded by $ 2 c_1 \big(\frac{\lambda}{a_{j}}\big)^4  $ and  $ 2 c_1 \big(\frac{\lambda}{a_{j+1}}\big)^4 $ respectively.

  Combining the above bounds for sufficiently large $j$ we get 
  \begin{align*} 
& \Big| J(a_{j})-  J(a_{j+1})  -\alpha c_{3} s_{0}^{-4\alpha-2}\cdot (a_{j}^{-2}-a_{j+1}^{-2}) \Big| \\
&\leq s_{0}^{-4\alpha}\int_{-A}^{A}|\widehat\psi(\lambda)|^{2} \Bigg[  \frac{2^7}{3^3}\Bigg(\frac{\lambda}{a_{j}s_{0}}\Bigg)^4+ \frac{2^7}{3^3} \Bigg(\frac{\lambda}{a_{j+1}s_{0}}\Bigg)^4 +   2 c_1\Bigg[  \Bigg(\frac{\lambda}{a_{j}}\Bigg)^4+  \Bigg(\frac{\lambda}{a_{j+1 }}\Bigg)^4\Bigg] \Bigg]d \lambda\\
& \leq \frac{2^8}{3^3}s_{0}^{-4\alpha-4} a_{j}^{-4} \int_{-A}^{A}|\widehat\psi(\lambda)|^{2}\lambda^4\bigg( 1+ 3^3 c_1 s_{0}^4/2^6  \bigg)d\lambda.
    \end{align*}  
Thus, for sufficiently large $j$ we obtain
 \begin{align*} 
 &\Bigg|\frac {J(a_{j})-  J(a_{j+1})} {a_{j}^{-2}-a_{j+1}^{-2} } -\alpha c_{3} s_{0}^{-4\alpha-2}  \Bigg|\leq\frac{ 2^{7}A^{2 }s_{0}^{-4\alpha-4} \left(1+3^3 c_1 s_{0}^4/2^6 \right) }{3^3 a^2_{j}\Big(1-\frac{a_{j}^2}{ a_{j+1}^2} \Big)}    c_{3}  
\end{align*} 
which completes the proof of Lemma  \ref{lemma_7}.\hfill\(\Box\)

Now let us define $\Delta \bar{\delta}_{ j\cdot}^{(2)} =  \frac{\bar{\delta}_{ j\cdot}^{(2)}-\bar{\delta}_{ j+1\cdot}^{(2)}} {a_{j}^{-2} -a_{j+1}^{-2}   }.$ Then the following result holds.
\begin{proposition}\label{pro_2}
Let the assumptions of Lemma {\rm \ref{lemma_3}} hold true and there exist $\varepsilon>0$ and $j_{0}\in \mathbb{N}$ such that $a_{j+1}\geq(1+\varepsilon) a_{j}$ for all $j\geq j_{0}.$  Then\\
 \[\Delta \bar{\delta}_{ j\cdot}^{(2)} \xrightarrow{a.s.} \alpha  c_{3} s_{0}^{-4\alpha-2}, \quad j \to+\infty.\]\\
Moreover, there exists an almost surely finite random variable $   c_{9}$ such that  for all $j\in \mathbb{N}$ it holds 
 \begin{align} \label{4_3}
\Big| \Delta \bar{\delta}_{ j\cdot}^{(2)} - \alpha  c_{3} s_{0}^{-4\alpha-2}  \Big| \leq   c_{9} \max\Big (a_{j}^{2} r_{j}, a_{j}^{-2}\Big). 
\end{align}  
\end{proposition}  
 \noindent\textit{Proof}
 Note that  $\Delta \bar{\delta}_{ j\cdot}^{(2)}$ can be rewritten as
   \[\Delta \bar{\delta}_{ j\cdot}^{(2)}= \frac{ \big(\bar{\delta}_{ j\cdot}^{(2)}-J(a_{j})\big)  -  \big(\bar{\delta}_{ j+1\cdot}^{(2)}-J(a_{j+1})\big)   + \big(  J(a_{j})- J(a_{j+1})\big) }    {a_{j}^{-2} -a_{j+1}^{-2}}.
 \]
 Thus, by Lemmas  \ref{lemma_3} and  \ref{lemma_7}
 \begin{align*} 
  &\Big| \Delta \bar{\delta}_{ j\cdot}^{(2)} - \alpha c_{3} s_{0}^{-4\alpha-2}   \Big| \leq \frac{c_{6}(r_{j}+r_{j+1}) }{a_{j}^{-2} -a_{j+1}^{-2}} + \frac{ 2^{6}A^{2 }s_{0}^{-4\alpha-4}\left(1+3^3 c_1 s^4_0/2^6\right)} {3^3 a^2_{j} \Big(1-\Big(\frac{a_{j}^2}{ a_{j+1}^2} \Big)^2\Big)}c_{3}\\
  &\leq \frac{2 c_{6}} {1-\big(\frac{a_{j}}  {a_{j+1}}\big)^2  }a^{2}_{j} r_{j}+ \frac{ 2^{6}A^{2 }s_{0}^{-4\alpha-4} \left(1+3^3c_1 s^4_0/2^6 \right)c_{3}}   {3^3\left(1-\big(\frac{a_{j}}  {a_{j+1}}\big)^2\right) } a_{j}^{-2}\\
  &\leq  \frac{2 c_{6}a^{2}_{j} r_{j}}   {    1-\frac{1}  {(1+\varepsilon)^2}}+ \frac{ 2^{6}A^{2 }s_{0}^{-4\alpha-4} \left(1+3^3 c_1 s^4_0/2^6\right)c_{3}}   {3^3\left(1-\frac{1}  {\left(1+\varepsilon\right)^2}\right) } a_{j}^{-2}\leq   c_{9}\max\Big (a_{j}^{2} r_{j}, a_{j}^{-2}\Big), 
    \end{align*} 
where $c_{9}:=\frac{2 c_{6}}   {    1-\frac{1}  {(1+\varepsilon)^2}}+ \frac{ 2^{6}A^{2 }s_{0}^{-4\alpha-4} \left(1+3^3 c_1 s^4_0/2^6\right)c_{3}}   {3^3\left(1-\frac{1}  {\left(1+\varepsilon\right)^2}\right) }. $ \hfill\(\Box\) 
\begin{remark}\label{remark_5}  
 As the rate of decay of $\{  r_{j}\}$ can be arbitrary selected the best upper bound given by $( \ref{4_3})$  has order $a^{-2}_{j}.$
\end{remark}
\section{Estimation of $(s_{0},\alpha)$}\label{sec_{5}}
In the previous sections we proved that if the true values of parameters are $(s_{0},\alpha)$, then the vector statistics
 \begin{align*}
\begin{pmatrix} 
\bar{\delta}_{ j\cdot}^{(2)}/   c_{2} \\
\Delta \bar{\delta}_{ j\cdot}^{(2)}/ c_{3}
\end{pmatrix}
 \xrightarrow{a.s.} 
\begin{pmatrix} 
s_{0}^{-4\alpha} \\
\alpha s_{0}^{-4\alpha-2}
\end{pmatrix}
,\quad j\to +\infty.  \end{align*}
 In this section we investigate properties of the pair $(\hat{s}_{0j}, \hat{\alpha }_{j})$ that is a solution of the system
   \begin{align}\label{5_1}
  \begin{cases}
        \hat{s}_{0j}^{-4 \hat{\alpha }_{j}}=\bar{\delta}_{ j\cdot}^{(2)}/   c_{2},\\
        \hat{\alpha }_{j}\hat{s}_{0j}^{-4\hat{\alpha }_{j}-2}=\Delta \bar{\delta}_{ j\cdot}^{(2)}/ c_{3}.
        \end{cases}
    \end{align}
 To handle the cases, where $\left(\frac{\bar{\delta}_{ j\cdot}^{(2)}}{   c_{2}}, \frac{\Delta \bar{\delta}_{ j\cdot}^{(2)}} { c_{3}}\right) $ may not be in the feasible region of $ \left( s_{0}^{-4 \alpha},   \alpha s_{0}^{-4\alpha -2} \right),$ we propose adjusted estimates.

First we discuss existence of solutions.
 \begin{lemma}\label{lemma_8} 
 Let $(y_1,y_2 )\in R_{\mathbf{y}},$ where $R_{\mathbf{y}}:=\Big\{ (y_1,y_2) \in (0,1)\times (0,y^2_{1}/2) \Big\}.$
Then the system 
    \begin{align} 
  \begin{cases}\label{5_2}
        s_{0}^{-4 \alpha}=y_1,\\
        \alpha s_{0}^{-4\alpha -2}=y_2,
        \end{cases}
    \end{align}
  has a solution $(s_{0}, \alpha)\in  (1,+\infty) \times \big(0,\frac{1}{2}\big)  .$ 
    \end{lemma}
  \noindent\textit{Proof}
  Let us find the range of the 2-d valued function 
\[ \mathbf{y}(s_{0},\alpha)=
\begin{pmatrix} 
     y_1(s_{0},\alpha)\\
        y_2(s_{0},\alpha)
\end{pmatrix}
 :=
\begin{pmatrix} 
s_{0}^{-4\alpha} \\
\alpha s_{0}^{-4\alpha-2}
\end{pmatrix}
\]
 defined on the domain  $(s_{0}, \alpha)\in (1,+\infty) \times  \big(0,\frac{1}{2}\big).$
 
 For simplicity, we use the notations $y_1$ and $y_2$ instead of $     y_1(s_{0},\alpha)$ and
  $  y_2(s_{0},\alpha)$ for the following computations.

  As $s_{0}>1,$ then for each $\alpha \in\big(0,\frac{1}{2}\big)$ the range of possible values of $y_1$ is $ (0,1).$  For each  $\alpha \in\big(0,\frac{1}{2}\big)$ and 
 $y_1\in (0,1)$ there is a such $s_{0}$ that $y_1=s_{0}^{-4 \alpha }$. The variable $y_2$ can be expressed in terms of $y_1$ as $y_2=\alpha \cdot y_1^{1+\frac{1}{2\alpha}}.$ Therefore, we can assume that $y_1$ is fixed and change only $\alpha$ to investigate the range of $y_2.$

 Notice, that 
 \[(y_2)'_{\alpha}=y_1^{1+\frac{1}{2\alpha}}+\alpha y_1^{1+\frac{1}{2\alpha}} \ln{(y_1)}\cdot \bigg(-\frac{1}{2\alpha^2}\bigg) =  y_1^{1+\frac{1}{2\alpha}}\bigg(1-\frac{\ln{(y_1)}}{2\alpha}\bigg).  \]

If $y_1\in (0,1),$ then $(y_2)'_{\alpha}>0$ and  $y_2$ is an increasing function of $\alpha$ with the range $ (0, \ y^2_{1}/2).$ 
Hence, the range of the function $\mathbf{y}(s_{0}, \alpha)$ on the domain $(1,+\infty)\times \big(0, \frac{1}{2}\big)$ is $ R_{\mathbf{y}}$, which completes the proof. \hfill \(\Box\)

 Thus, if $\bigg(\frac{\bar{\delta}_{ j\cdot}^{(2)}}{   c_{2}},\frac{\Delta \bar{\delta}_{ j\cdot}^{(2)}}{ c_{3}}\bigg)\in R_{\mathbf{y}}$ then there is a pair $(\hat{s}_{0j}, \hat{\alpha }_{j})\in (1,+\infty)\times(0, \frac{1}{2})$ that satisfies the system of equations $( \ref{5_1}).$
 Now we will investigate uniqueness of solutions.
 \begin{lemma}
 Let $(y_1,y_2 )\in R_{\mathbf{y}}.$ Then system $( \ref{5_2})$ has a unique solution.
 \end {lemma}
\noindent\textit{Proof} If $( s_{0},\alpha)$ and $(s'_{0},\alpha')$ are two solutions of the system $( \ref{5_2})$ for some $ (y_1,y_2 )\in R_{\mathbf{y}}$ then 
 \[ 
  \begin{cases}
        s_{0}^{-4 \alpha }= (s'_{0})^{-4 \alpha' },\\ 
        \alpha s_{0}^{-2}= \alpha' (s'_{0})^{-2},
 \end{cases}
\]
and therefore
\[
\begin{cases}
        s_{0}^{ \alpha }=  (s'_{0})^{ \alpha' },\\ 
       \Big(\frac{\alpha} {s_{0}^{ 2}}\Big)^{ \alpha }=  \Big(\frac {\alpha'} {(s'_{0})^{2}}\Big)^{ \alpha }.
 \end{cases}
 \]
 
Hence, $\Big(\frac {\alpha'} {\alpha }\Big)^{ \alpha }=\Big(\frac{(s'_{0})^{2}} {s_{0}^2}\Big)^{ \alpha }= (s'_{0})^{2(  \alpha-  \alpha' )}$ and $ \frac {\alpha'} {\alpha }=  (s'_{0})^{2(1- \frac {\alpha'} {\alpha })}.$

Denoting $\frac {\alpha'} {\alpha }=t,$ and $(s'_{0})^{2}=a$ we obtain the equation
    \begin{align}\label{5_3}
     t a ^{t-1}=1,\quad t \in \mathbb{R_+}.    \end{align} 
     
  As $s_0>1$ then $s'_0$ must also be greater than 1 (otherwise $\alpha'<0,$ which is not feasible). Hence,  $a >1$ and the left-hand side of (\ref{5_3}) is an increasing function. Hence, the equation has the  only solution $t=1,$ which means $\alpha'=\alpha$ and implies a unique solution of $( \ref{5_2}) $.  \hfill\(\Box\)

Now we provide solutions to system (\ref{5_2}). These solutions are given in terms of the $LambertW$ function, which is defined as a solution of the equation  
\[t \mathrm{e}^{t}=x, \quad x\geq \mathrm{e}^{-1}, \]
 i.e. $t = LambertW(x) .$  

 \begin{proposition}\label{pro_3}
 Let $(y_1,y_2)\in R_{\mathbf{y}}.$ Then the solution to system $( \ref{5_2})$ is \begin{align}\label{5_4} 
\begin{array}{ll}
s_{0}&=exp\Big (\frac{1}{2}  LambertW\Big(\frac{y_1} {y_2}\ln\Big(y_1^{-\frac{1}{2}}\Big) \Big)    \Big),\\
\alpha&= \frac{y_2} {y_1} exp\Big (  LambertW\Big(\frac{y_1} {y_2}\ln\Big(y_1^{-\frac{1}{2}}\Big) \Big) \Big).
\end{array}
\end{align} 
 \end{proposition}
\noindent\textit{Proof} Let us rewrite (\ref{5_2}) as 
 \[
 \begin{cases}
       -4 \alpha\ln( s_{0})=  \ln(y_1),\\ 
      \alpha=\frac{y_2} {y_1} s^2_{0},  
 \end{cases}
\]
and therefore
\[
 \begin{cases}
       \alpha\big(\ln( \alpha)+\ln\big(\frac{y_1} {y_2}\big)\big)= -\frac{ \ln (y_1)}{2},\\ 
      s_{0} = \sqrt{\frac{y_1} {y_2}\alpha}. 
 \end{cases}
\]    
 Denoting $t=\ln(\alpha)$ the first equation can be rewritten as

 \[ \mathrm{e}^{t}\left(t+\ln\left(\frac{y_1} {y_2}\right)\right)=- \frac{ \ln (y_1)}{2},\]

 \[ \mathrm{e}^{t+\ln\Big(\frac{y_1} {y_2}\Big)} \Big(t+\ln\Big(\frac{y_1} {y_2}\Big)\Big)= \mathrm{e}^{\ln\Big(\frac{y_1} {y_2}\Big)} \ln \left(y^{-\frac{1}{2}}_{1}\right)= \frac{y_1}{y_2} \ln \left(y^{-\frac{1}{2}}_{1}\right).\]
 Hence, by the definition of the LambertW function we obtain
 \[  
t= LambertW\Big(\frac{y_1} {y_2}\ln \Big(y_1^{-\frac{1}{2}}\Big)\Big)-\ln\Big(\frac{y_1} {y_2}\Big).
\]
Finally, $( \ref{5_4})$ follows from $\alpha=  \mathrm{e}^{t}$ and $ s_{0}=\Big(\frac{y_1}{y_2}\alpha\Big)^{\frac{1}{2}}.  $ \hfill \(\Box\)
 \begin{remark}\label{remark_8}
 The LambertW function has two real branches Lambert$W_{0}$ and Lambert$W_{-1}.$ The branch Lambert$W_0$ is defined on the interval $[-\frac{1}{\mathrm{e}},+\infty),$ but the branch Lambert$W_{-1}$ is defined only on  the interval $ \left[-\frac{1}{\mathrm{e}},0\right).$ The point $ (-\frac{1}{\mathrm{e}},-1) $ is a branch point for Lambert$W_{0}$  and Lambert$W_{-1}.$ 
 
Hence, for  $y_{1}\in(0,1)$ 
it holds  $ \Big(\frac{y_1} {y_2}\Big)\ln\Big(y_1^{-\frac{1}{2}}\Big)>0$ and $ ( \ref{5_4}) $ gives a unique solution to {\rm(\ref{5_2})} with the branch Lambert$W_0.$
  \end{remark}

 Now, we see that for $(y_1,y_2) \in R_{\mathbf{y}}$ there is a unique solution $(s_{0},{\alpha})$ to (\ref{5_2}).
If $ s_{0}$ and $\alpha$ are the true value of parameters then the corresponding $(y_1,y_2)\in R_{\mathbf{y}}.$ As $R_{\mathbf{y}}$ is an open set, then  $(y_1,y_2)\in int(R_{\mathbf{y}})=R_{\mathbf{y}}$ and there is some $j_{0}\in \mathbb{N}$ such that $\bigg(\frac{\bar{\delta}_{ j\cdot}^{(2)}}{   c_{2}},\frac{\Delta \bar{\delta}_{ j\cdot}^{(2)}}{ c_{3}}\bigg)\in R_{\mathbf{y}}$  for all $j\geq j_{0\cdot},$ where   $int(\cdot)$ denotes the interior of a set. Therefore, starting from $  j_{0}$ system (\ref{5_1}) has a unique solution.

However, it might happen that $\bigg(\frac{\bar{\delta}_{ j\cdot}^{(2)}}{   c_{2}},\frac{\Delta \bar{\delta}_{ j\cdot}^{(2)}}{ c_{3}}\bigg)\notin R_{\mathbf{y}}$ for some $j<j_0$ even if   $(y_1,y_2)\in R_{\mathbf{y}}=int( R_{\mathbf{y}})$ for the  corresponding true value $ (s_{0},\alpha).$   For the cases  $\bigg(\frac{\bar{\delta}_{ j\cdot}^{(2)}}{   c_{2}},\frac{\Delta \bar{\delta}_{ j\cdot}^{(2)}}{ c_{3}}\bigg) \notin R_{\mathbf{y}}$  to define $(\hat{s}_{0j}, \hat{\alpha}_j )$ we introduce "adjusted" values $\bigg(\frac{\bar{\delta}_{ j\cdot}^{(2,a)}}{   c_{2}},\frac{\Delta \bar{\delta}_{ j\cdot}^{(2,a)}}{ c_{3}}\bigg) \in  R_{\mathbf{y}}.$  
\begin{Definition}\label{5_1}
The adjusted statistics
 $ \frac{\bar{\delta}_{ j\cdot}^{(2,a)}}{   c_{2}}$ and $\frac{\Delta \bar{\delta}_{ j\cdot}^{(2,a)}}{ c_{3}}  $ are  defined as follows
 \begin{itemize}
 \item if $\frac{\bar{\delta}_{ j\cdot}^{(2)}}{c_{2}} \in(0,1)$ and $\frac{\Delta \bar{\delta}_{ j\cdot}^{(2)}}{ c_{3}} \ge \frac{1}{2}  \Big(\frac{\bar{\delta}_{ j\cdot}^{(2)}}{   c_{2}}\Big)^2,$ then
  \[\bar{\delta}_{ j\cdot}^{(2,a)}= \bar{\delta}_{ j\cdot}^{(2)}\quad \text{and} \quad \Delta \bar{\delta}_{ j\cdot}^{(2,a)}= c_{3}\max \left( \left(\frac{\bar{\delta}_{ j\cdot}^{(2)}}{   c_{2}}\right)^2 - \frac{\Delta \bar{\delta}_{ j\cdot}^{(2)}}{ c_{3}}, \frac{1}{4}\left(\frac{\bar{\delta}_{ j\cdot}^{(2)}}{   c_{2}}\right)^2 \right);\] 	

 \item if $\frac{\bar{\delta}_{ j\cdot}^{(2)}}{   c_{2}} \in(0,1)$ and $\frac{\Delta \bar{\delta}_{ j\cdot}^{(2)}}{ c_{3}} \le 0,$ then\\ \[\bar{\delta}_{ j\cdot}^{(2,a)}= \bar{\delta}_{ j\cdot}^{(2)}\quad \text{and} \quad \Delta \bar{\delta}_{ j\cdot}^{(2,a)}=c_{3}\min \left( - \frac{\Delta \bar{\delta}_{ j\cdot}^{(2)}}{ c_{3}}, \frac{1}{4}\left(\frac{\bar{\delta}_{ j\cdot}^{(2)}}{   c_{2}}\right)^2 \right);\]  
 
  \item if $\frac{\bar{\delta}_{ j\cdot}^{(2)}}{c_{2}} \ge 1$ and $\frac{\Delta \bar{\delta}_{ j\cdot}^{(2)}}{ c_{3}} <\frac{1}{2},$ then\\ \[\Delta \bar{\delta}_{ j\cdot}^{(2,a)}=\Delta \bar{\delta}_{ j\cdot}^{(2)}\quad \text{and} \quad  \bar{\delta}_{ j\cdot}^{(2,a)}= c_2\max\left(2-\frac{\bar{\delta}_{ j\cdot}^{(2)}}{c_{2}},\frac{1}{2}\left(1+c_2\left(2\frac{\Delta \bar{\delta}_{ j\cdot}^{(2)}}{c_3}\right)^\frac{1}{2}\right)\right);\] 
  
    \item if $\frac{\bar{\delta}_{ j\cdot}^{(2)}}{c_{2}} \ge 1$ and $\frac{\Delta \bar{\delta}_{ j\cdot}^{(2)}}{ c_{3}} \ge\frac{1}{2},$ then\\ \[\Delta \bar{\delta}_{ j\cdot}^{(2,a)}=c_3\max\left(1-\frac{\Delta \bar{\delta}_{ j\cdot}^{(2)}}{ c_{3}},\frac{1}{4}\right)\quad \text{and} \quad  \bar{\delta}_{ j\cdot}^{(2,a)}= c_2\max\left(2-\frac{\bar{\delta}_{ j\cdot}^{(2)}}{c_{2}},\frac{1}{2}\left(1+c_2\left(2\frac{\Delta \bar{\delta}_{ j\cdot}^{(2,a)}}{c_3}\right)^\frac{1}{2}\right)\right);\]   	
 	
    \item if $\frac{\bar{\delta}_{ j\cdot}^{(2)}}{c_{2}} \ge 1$ and $\frac{\Delta \bar{\delta}_{ j\cdot}^{(2)}}{ c_{3}} \le 0,$ then\\ \[\Delta \bar{\delta}_{ j\cdot}^{(2,a)}=c_3\min\left(-\frac{\Delta \bar{\delta}_{ j\cdot}^{(2)}}{ c_{3}},\frac{1}{4}\right)\quad \text{and} \quad  \bar{\delta}_{ j\cdot}^{(2,a)}= c_2\max\left(2-\frac{\bar{\delta}_{ j\cdot}^{(2)}}{c_{2}},\frac{1}{2}\left(1+c_2\left(2\frac{\Delta \bar{\delta}_{ j\cdot}^{(2,a)}}{c_3}\right)^\frac{1}{2}\right)\right);\] 
     
\item otherwise $\bar{\delta}_{ j\cdot}^{(2,a)}=\bar{\delta}_{ j\cdot}^{(2)}$ and $\Delta \bar{\delta}_{ j\cdot}^{(2,a)}= \Delta \bar{\delta}_{ j\cdot}^{(2)}.$
  \end{itemize}
\end{Definition}
In the fourth and fifth cases the value of $\Delta \bar{\delta}_{j\cdot}^{(2,a)}$ is computed first and then it is used to compute the adjusted value $\bar{\delta}_{ j\cdot}^{(2,a)}.$

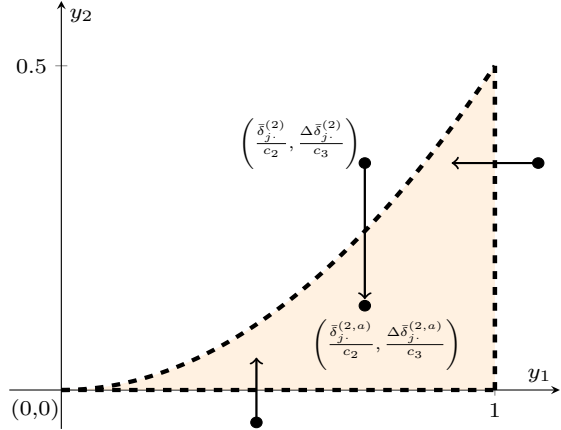
\begin{figure}[h] 
	\begin{center}
\resizebox{0.5\textwidth}{1\height}{	
	\begin{tikzpicture}
	\begin{axis}[
	axis lines = middle,
	xlabel = $y_1$,
	ylabel = {$y_2$},xmin=-0.12,
	xmax=1.15,ymin=-0.06,
	ymax=0.6, xtick={0,1},
	ytick={0,0.5}
	]
	
	\addplot [name path=axis,
	domain=0:1, 
	samples=100, 
	dashed,style={ultra thick}
	]
	{0};
	
	\addplot [name path=f,
	domain=0:1, 
	samples=100, 
	dashed,style={ultra thick}
	]
	{0.5*x^2};
	\node at (axis cs:  -0.06, -0.03)
	{(0,0)};
	\node at (axis cs:  .55,  .38) {\scalebox{0.8}{$\left(\frac{\bar{\delta}_{ j\cdot}^{(2)}}{   c_{2}}, \frac{\Delta \bar{\delta}_{ j\cdot}^{(2)}}{ c_{3}}\right)$}};
	\node at (axis cs:  .75,  .07) {\scalebox{0.8}{$\left(\frac{\bar{\delta}_{ j\cdot}^{(2,a)}}{   c_{2}}, \frac{\Delta \bar{\delta}_{ j\cdot}^{(2,a)}}{ c_{3}}\right)$}};
	
	\addplot[dashed, style={ultra thick}]
	coordinates{(1,0) (1,0.5)};
	
	\addplot [
	thick,
	color=black,
	fill=bisque, 
	fill opacity=0.5
	]
	fill between[
	of=f and axis,
	soft clip={domain=0:1},
	];
	
	\addplot[mark=*] coordinates {(1.1,0.35)};
	\addplot[->,style={thick}] coordinates
	{(1.1,0.35) (0.9,0.35)};
	
	\addplot[mark=*] coordinates {(0.45,-0.05)};
	\addplot[->,style={thick}] coordinates
	{(0.45,-0.05) (0.45,0.05)};
	
	\addplot[mark=*] coordinates {(0.7,0.35)};
	\addplot[mark=*] coordinates {(0.7,0.13)};
	\addplot[->,style={thick}] coordinates
	{(0.7,0.35) (0.7,0.14)};
	\end{axis}
	\end{tikzpicture}
}
		\caption{Plot of $R_{\mathbf{y}},$ $\left(\frac{\bar{\delta}_{ j\cdot}^{(2)}}{   c_{2}}, \frac{\Delta \bar{\delta}_{ j\cdot}^{(2)}}{ c_{3}}\right)$ and  the corresponding $\left(\frac{\bar{\delta}_{ j\cdot}^{(2,a)}}{   c_{2}}, \frac{\Delta \bar{\delta}_{ j\cdot}^{(2,a)}}{ c_{3}}\right).$ } \label{fig2}
\end{center}
\end{figure}

Figure~\ref{fig2} clarifies geometric reasons to introduce the adjusted values $ \left(\frac{\bar{\delta}_{ j\cdot}^{(2,a)}}{   c_{2}}, \frac{\Delta \bar{\delta}_{ j\cdot}^{(2,a)}}{ c_{3}}\right).$ Vertical or horizontal reflections over  boundaries of  $R_{\mathbf{y}}$ are used with an additional constraint that the reflected points do not go beyond  the opposite boundaries of  $R_{\mathbf{y}}.$ Also, the reflected points should not belong to the boundaries.  For instance, this might happen, in the first case of Definition \ref{5_1} if one has $ \frac{\Delta \bar{\delta}_{ j\cdot}^{(2)}}{ c_{3}}  = \frac{1}{2} \left(\frac{\bar{\delta}_{ j\cdot}^{(2)}}{   c_{2}}\right)^2 .$  

\begin{remark}\label{remark_9}
By the construction in Definition \rm \ref{5_1} the adjusted pair $ \left(\frac{\bar{\delta}_{ j\cdot}^{(2,a)}}{   c_{2}}, \frac{\Delta \bar{\delta}_{ j\cdot}^{(2,a)}}{ c_{3}}\right) \in R_{\mathbf{y}}$  and the both $\bigg(\frac{\bar{\delta}_{ j\cdot}^{(2,a)}}{   c_{2}},\frac{\Delta \bar{\delta}_{ j\cdot}^{(2,a)}}{ c_{3}}\bigg)$ and $\bigg(\frac{\bar{\delta}_{ j\cdot}^{(2)}}{   c_{2}},\frac{\Delta \bar{\delta}_{ j\cdot}^{(2)}}{ c_{3}}\bigg)$ converge to the same value  $(s_{0}^{-4\alpha}, \alpha s_{0}^{-4\alpha-2})$  when $j\to+\infty.$ 
\end{remark}

\begin{remark}\label{remark_10}
As only a finite number of $\left(\frac{\bar{\delta}_{ j\cdot}^{(2)}}{   c_{2}}, \frac{\Delta \bar{\delta}_{ j\cdot}^{(2)}}{ c_{3}}\right)$ fall outside of $ R_{\mathbf{y}},$ then there is $j_0 \in \mathbb{N}$ such that $ \left(\frac{\bar{\delta}_{ j\cdot}^{(2,a)}}{   c_{2}}, \frac{\Delta \bar{\delta}_{ j\cdot}^{(2,a)}}{ c_{3}}\right) = \left(\frac{\bar{\delta}_{ j\cdot}^{(2)}}{   c_{2}}, \frac{\Delta \bar{\delta}_{ j\cdot}^{(2)}}{ c_{3}}\right) $ for all $j\geq j_0.$
Therefore, in this case $ \left(\frac{\bar{\delta}_{ j\cdot}^{(2,a)}}{   c_{2}}, \frac{\Delta \bar{\delta}_{ j\cdot}^{(2,a)}}{ c_{3}}\right)$ and $ \left(\frac{\bar{\delta}_{ j\cdot}^{(2)}}{   c_{2}}, \frac{\Delta \bar{\delta}_{ j\cdot}^{(2)}}{ c_{3}}\right)$ have the same rate of convergence to $(s_{0}^{-4\alpha}, \alpha s_{0}^{-4\alpha-2})$ when $j\to +\infty.$
\end{remark}

Now we are ready to formulate the main result.
  \begin{theorem}\label{Theorem_1} 
 Let the process $X(t)$ and the filter $\psi(\cdot)$ satisfy Assumptions \rm{\ref{Assumption_1}--\ref{Assumption_3}}.
 \it{Let $(\hat{s}_{0j}, \hat{\alpha }_{j})$ be a solution of the system of equations
 \begin{align*} 
\begin{cases}
   \hat{\alpha }_{j}(\hat{s}_{0j})^{-4\hat{\alpha }_{j}-2}=\Delta \bar{\delta}_{ j\cdot}^{(2,a)}/ c_{3},\\
  \hat{s}_{0j}^{-4 \hat{\alpha }_{j}}=\bar{\delta}_{ j\cdot}^{(2,a)}/   c_{2},\\
     \end{cases}        
\end{align*}
 where  $\Delta \bar{\delta}_{ j\cdot}^{(2,a)}$  and $ \bar{\delta}_{ j\cdot}^{(2,a)}$ are the adjusted statistics.  
 Then
 \begin{align}\label{5_6} 
\begin{array}{ll}
 \widehat s_{0j} &=exp\left (\frac{1}{2} LambertW\left(\frac{\ln\big(   c_{2}/ \bar{\delta}_{ j\cdot}^{(2,a)}   \big)} {2q_j} \right)  \right) , \\
\widehat \alpha_j &= q_j\ exp \left(  LambertW\left(\frac{\ln\big(   c_{2}/ \bar{\delta}_{ j\cdot}^{(2,a)}   \big)} {2q_j}\right)  \right), 
\end{array}
\end{align}
 where $q_j= \frac{   c_{2}}{ c_{3}}  \frac{\Delta \bar{\delta}_{ j\cdot}^{(2,a)}}{ \bar{\delta}_{ j\cdot}^{(2,a)}}.$

If $s_{0}$ and $\alpha$ are the true values of parameters and the assumptions of Proposition \rm{\ref  {pro_2}} \textit{ hold true, then
 $\widehat s_{0j} \xrightarrow{a.s.} s_{0}$ and $ \widehat \alpha_j \xrightarrow{a.s.}\alpha,$ when $ j\to +\infty.$
 Moreover, there are almost surely finite random variables $c_{10}$ and $c_{11}$ such that for all $j\in \mathbb{N}$ it holds
 \[
      |\widehat s_{0j}-s_{0}|\leq c_{10} \max\big (a_{j}^{2} r_{j}, a_{j}^{-1}\big)\
      \] \text{and} \[|\widehat \alpha_{j}-\alpha|\leq c_{11} \max\big (a_{j}^{2} r_{j}, a_{j}^{-1}\big). 
      \]}}
  \end{theorem}  
  \noindent\textit{Proof}
The first statement of the theorem immediately follows from Proposition \ref{pro_3}.

To investigate properties of the solutions $(\hat{s}_{0j}, \hat{\alpha }_{j})$ one has to study properties of $LambertW\big(\frac{x_{1}}{x_{2}}\big).$

Notice that, by the 2-dimensional mean value theorem, for any $f:\mathbb{R }^{2}\to \mathbb{R}$ that is from $C^1(\mathbb{R}^2)$ it holds
      \begin{align*}
      f(\mathbf{x})-f(\mathbf{y})=\left(\nabla f\left((1-c)\mathbf{x}+c\mathbf{y}\right), \mathbf{x}-\mathbf{y}\right), \quad\mathbf{ x,y}\in \mathbb{R}^{2}, 
        \end{align*}
where $c\in[0,1], \nabla$ denotes the gradient and $(\cdot,\cdot)$ is the scalar product in $\mathbb{R}^{2}.$
Therefore, 
        \begin{align*}
             \big |f(\mathbf{x})-f(\mathbf{y})\big|\leq \sup_{c\in[0,1]}\big|\big|\nabla f\big((1-c)\mathbf{x}+c\mathbf{y}\big|\big|\cdot\big|\big|\mathbf{x}-\mathbf{y}\big|\big|.
          \end{align*}
Now, applying this result to the function $f(\mathbf{x})=\exp\left(\frac{1}{2}  LambertW\big(\frac{x_{1}}{x_{2}}\big) \right),$  $\mathbf{x}=~\Big(\ln\Big(\frac{   c_{2}}{ \bar{\delta}_{ j\cdot}^{(2,a)}}\Big), 2q_j \Big),$ \\ $\mathbf{y}=\left(4\alpha\ln(s_{0}), \frac{2\alpha}{s_{0}^2}\right),$ and noting that the solution $\widehat s_{0j}$ is given by (\ref{5_6}), we obtain         
  \begin{align}\label{5_7} 
  |\widehat s_{0j}-s_{0}|&= \Big| \mathrm{e}^{\frac{1}{2} LambertW\big(\frac{x_{1}}{x_{2}}\big)}-  \mathrm{e}^{\frac{1}{2} LambertW\big(\frac{y_{1}}{y_{2}}\big)}\Big|
  \leq \sup_{c\in[0,1]}\Big|\Big|\nabla \mathrm{e}^{\frac{1}{2} LambertW\big(\frac{(1-c)x_{1}+cy_1}{(1-c)x_{2}+cy_2}\big)}\Big|\Big|\\
 &\times \Bigg|\Bigg|\Bigg(\ln\left(\frac{   c_{2}}{ \bar{\delta}_{ j\cdot}^{(2,a)}}\right) -4\alpha\ln(s_{0}),\frac{2   c_{2}}{ c_{3}}\frac{\Delta \bar{\delta}_{ j\cdot}^{(2,a)}}{\bar{\delta}_{ j\cdot}^{(2,a)}}-\frac{2\alpha}{s_{0}^2}\Bigg)\Bigg|\Bigg|.
  \end{align}
  Noting that $(LambertW(x))'=\frac{LambertW(x)}{x(1+LambertW(x))}$ we obtain
     \begin{align*}
     \nabla  \mathrm{e}^{\frac{1}{2} LambertW\big(\frac{x_{1}}{x_{2}}\big)}&= \mathrm{e}^{\frac{1}{2} LambertW\big(\frac{x_{1}}{x_{2}}\big)}\frac { LambertW\big(\frac{x_{1}}{x_{2}}\big)}{2\frac{x_1}{x_2}\big(1+LamberW\big(\frac{x_{1}}{x_{2}}\big)\big)} \Big(\frac{1}{x_2},-\frac{x_1}{x_2^2}  \Big)\\
     &= \mathrm{e}^{\frac{1}{2} LambertW\big(\frac{x_{1}}{x_{2}}\big)}\frac {LambertW\big(\frac{x_{1}}{x_{2}}\big)}{1+LambertW\big(\frac{x_{1}}{x_{2}}\big)} \Big(\frac{1}{2x_1},-\frac{1}{2x_2}  \Big).
       \end{align*}
   By the properties of the adjusted estimates  $\mathbf{x}= \Big(\ln\Big(\frac{   c_{2}}{\bar{\delta}_{ j\cdot}^{(2,a)}}\Big), 2q_j\Big)\to \Big (4\alpha \ln(s_{0}), \frac{2\alpha}{s_0}^2\Big),$ when $j \to \infty.$ As the both $4\alpha\ln(s_{0})$ and $\frac{2\alpha^2}{s_{0}}$ are different from zero, then $x_1$ and $x_2$ are non-zero values separated  from zero for sufficiently large $j.$ Hence, $\frac{1}{x_1}$ and $\frac{1}{x_1}$ are bounded and the above gradient is bounded.

Now, we study the second multiplier in (\ref{5_7})
    \begin{align}\label{5_8} 
   \Bigg|\Bigg|\Bigg(\ln\Bigg(\frac{   c_{2}}{ \bar{\delta}_{ j\cdot}^{(2,a)}}\Bigg) -4\alpha\ln(s_{0}),\frac{2   c_{2}}{ c_{3}}\frac{\Delta \bar{\delta}_{ j\cdot}^{(2,a)}}{\bar{\delta}_{ j\cdot}^{(2,a)}}-\frac{2\alpha}{s_{0}^2}\Bigg)\Bigg|\Bigg|
   &\leq\Bigg|\ln\Bigg(\frac{   c_{2}}{ \bar{\delta}_{ j\cdot}^{(2,a)}s_{0}^{4\alpha}}\Bigg)\Bigg|  \\+\Bigg|  \frac{2   c_{2} s_{0}^{2} \Delta \bar{\delta}_{ j\cdot}^{(2,a)} -2\alpha  c_{3} \bar{\delta}_{ j\cdot}^{(2,a)}} { c_{3} s_{0}^{2} \bar{\delta}_{ j\cdot}^{(2,a) } }\Bigg|.
 \end{align}
 As $|\ln(x_1)-\ln(y_1)|\leq\frac{|x_1-y_1|}{\min(x_1,y_1)} $ for $x_1, y_1 \in \mathbb{R_+},$ then by  Proposition \ref{pro_1} and Remark \ref{remark_10} we can estimate the first summand in  (\ref{5_8}) as 
\begin{align}\label{5_9}   
 \Bigg|\ln\Bigg(\frac{   c_{2}s_{0}^{-4\alpha}}{ \bar{\delta}_{ j\cdot}^{(2,a)}}\Bigg)\Bigg|\leq\frac{\big|\bar{\delta}_{ j\cdot}^{(2,a)}-     c_{2} s_{0}^{-4\alpha}\big|}{\min \big(\bar{\delta}_{ j\cdot}^{(2,a)},    c_{2} s_{0}^{-4\alpha} \big)}\leq c_{12}\max \big(r_j, a_j^{-1}\big),
   \end{align}
where $c_{12}$ is an almost surely finite random variable.

The second summand in (\ref{5_8})  
can be estimated using Propositions  \ref{pro_1}, \ref{pro_2}  and Remark \ref{remark_10} as 
      \begin{align}\label{5_10} 
   &\Bigg| \frac{2   c_{2} s_{0}^{2} \Delta \bar{\delta}_{ j\cdot}^{(2,a)} -2\alpha  c_{3} \bar{\delta}_{ j\cdot}^{(2,a)}} { c_{3} s_{0}^{2} \bar{\delta}_{ j\cdot}^{(2,a) } }   \Bigg|
 \leq    \Bigg| \frac{2   c_{2} s_{0}^{2}\big( \Delta \bar{\delta}_{ j\cdot}^{(2,a)} -2\alpha s_{0}^{-4\alpha-2}  c_{3}\big) } { c_{3} s_{0}^{2} \bar{\delta}_{ j\cdot}^{(2,a) } }   \Bigg| + \Bigg| \frac{2 \alpha  c_{3} \big( s_{0}^{-4\alpha-2}- \bar{\delta}_{ j\cdot}^{(2,a)}\big) }{ c_{3} s_0^2\bar{\delta}_{ j\cdot}^{(2,a) } }   \Bigg|\\
 & \leq c_{13}\big(\max \big(a^{2}_{j} r_j, a^{-2}_{j}\big)+\max \big(r_j, a_j^{-1}\big)\big) \leq 2 c_{13}\max \big(a^{2}_{j} r_j, a^{-1}_{j}\big),
      \end{align}
  where $ c_{13}$  is an almost surely finite random variable. 

 Combining the upper bounds in (\ref{5_9}) and (\ref{5_10}) and denoting $c_{10}=\max \big(c_{12}, 2 c_{13}\big)$
   we obtain $|\widehat s_{0j}-s_{0}|\leq c_{10} \max\big (a_{j}^{2} r_{j}, a_{j}^{-1}\big).$

Finally, noting that 
       \begin{align*} 
   |\widehat \alpha_{j}-\alpha|&= \bigg|\frac{x_2}{2}   \mathrm{e}^{ LambertW\big(\frac{x_{1}}{x_{2}}\big)}- \frac{y_2}{2}   \mathrm{e}^{ LambertW\big(\frac{y_{1}}{y_{2}}\big)}  \bigg |
  \\& \leq  \frac{1}{2} \bigg | \mathrm{e}^{ LambertW\big(\frac{x_{1}}{x_{2}}\big)} (x_2-y_2) \bigg |+ \frac{1}{2} \bigg |y_2\left( \mathrm{e}^{ LambertW\big(\frac{x_{1}}{x_{2}}\big)}-  \mathrm{e}^{ LambertW\big(\frac{y_{1}}{y_{2}}\big)}\right)  \bigg |
        \end{align*}
 and using the upper bounds in (\ref{5_9}) and (\ref{5_10}) we obtain $|\widehat \alpha_{j}-\alpha|\leq c_{11} \max\big (a_{j}^{2} r_{j}, a_{j}^{-1}\big),$ which completes the proof.\hfill \(\Box\)
 
\section{Simulation studies}\label{sec_6}

This section  presents some numerical studies to confirm the theoretical findings. The results demonstrate that the approach can be extended to other  processes and filters.  

All  theoretical results in the paper were developed for functional data. However, only discrete time can be used for computer simulations. Therefore, we selected a large simulation grid which makes simulation  results very close to the case of continuous-time realisations.

We consider the Gegenbauer random process $X(t), t\in \mathbb{Z},$ see \cite{Alomari:2017}, \cite{Espejo:2015}
  and the references therein. This random process satisfies the following  equation  
\[
 \Delta^{d} _{u}X(t):= (1-2uB+B^2)^{d} X(t)=\varepsilon_t,
\]
where $ \Delta^{d} _{u}$ is the fractional difference operator given by \[\Delta^{d} _{u}= (1-2uB+B^2)^{d}= (1-2\cos(\nu) B+B^2)^{d}=[(1-e^{i\nu}B) (1-e^{-i\nu}B)],\]
$B$ denotes the backward-shift operator for the time coordinate $t,$ i.e. $B X_t= X_{t-1},$ $u= \cos\nu$ (i.e. $\nu= \arccos (u), |u|\leq1),$  $d \in (-\frac{1}{2}, \frac{1}{2} ),$  and
$\varepsilon_t$ is a zero-mean white noise with the common variance $E (\varepsilon^2_t)= \sigma^2_{\varepsilon}.$

There exists the following representation of a stationary Gegenbauer random process
\begin{equation}
 X(t)= \sum_{n=0}^{\infty}  C_{n}^{(d)} (u) \varepsilon_{t-n },\label{GRF}\end{equation}
  where $d$ $\neq 0$  and the Gegenbauer polynomial $ C_{n}^{(d)}
(u)$ is given by
\begin{equation*}
 C_{n}^{(d)} (u)= \sum_{k=0}^{[n/2]} (-1)^{k}\frac{(2u)^{n-2k}\Gamma (d-k+n)}{k!(n-2k)!\Gamma (d)},
\end{equation*} 
where $[n/2]$ is the integer part of $n/2,$ and  $\Gamma(\cdot)$ is the gamma function.

We generated the random process $X(t)$ using the parameter values $d=0.1$ and $u=0.3.$ The chosen parameters $d $ and $u$ correspond to $\ s_{0}$ and $\alpha$ inside of the admissible region $R_\mathbf{y}.$  The realisations of $X(t)$  were approximated by truncated sums with 40 terms in (\ref{GRF}). The filter transform  of $X(t)$ defined by (\ref{2_1}) was computed using the R package {\sc wmtsa.}

\begin{figure}[!h] 
 	\begin{center}
 		\begin{minipage}{5.4cm} 
 			\includegraphics[width=\textwidth]{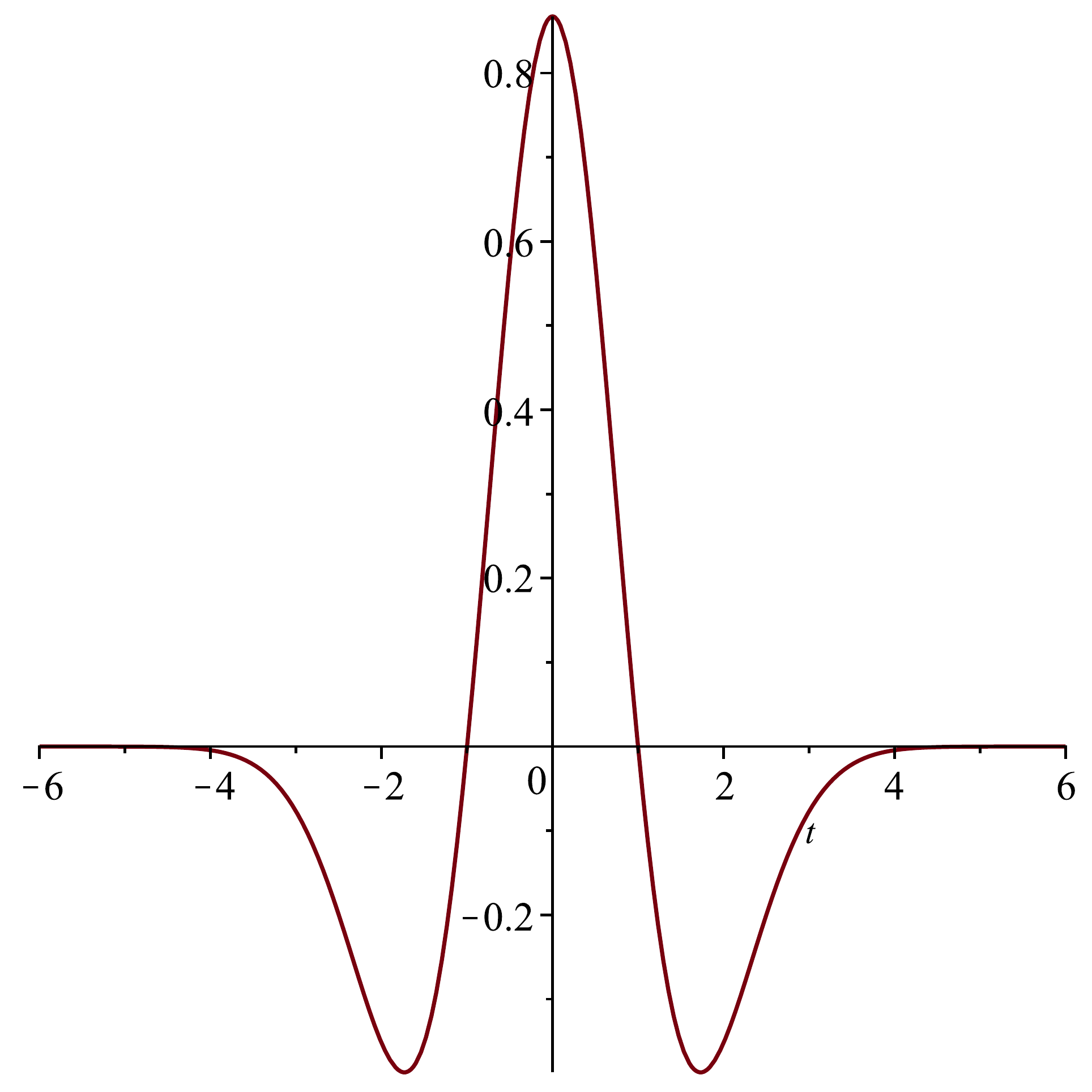}
 		\end{minipage}%
 		\quad \quad
 		\begin{minipage}{5.4cm}
 			\includegraphics[width=\textwidth]{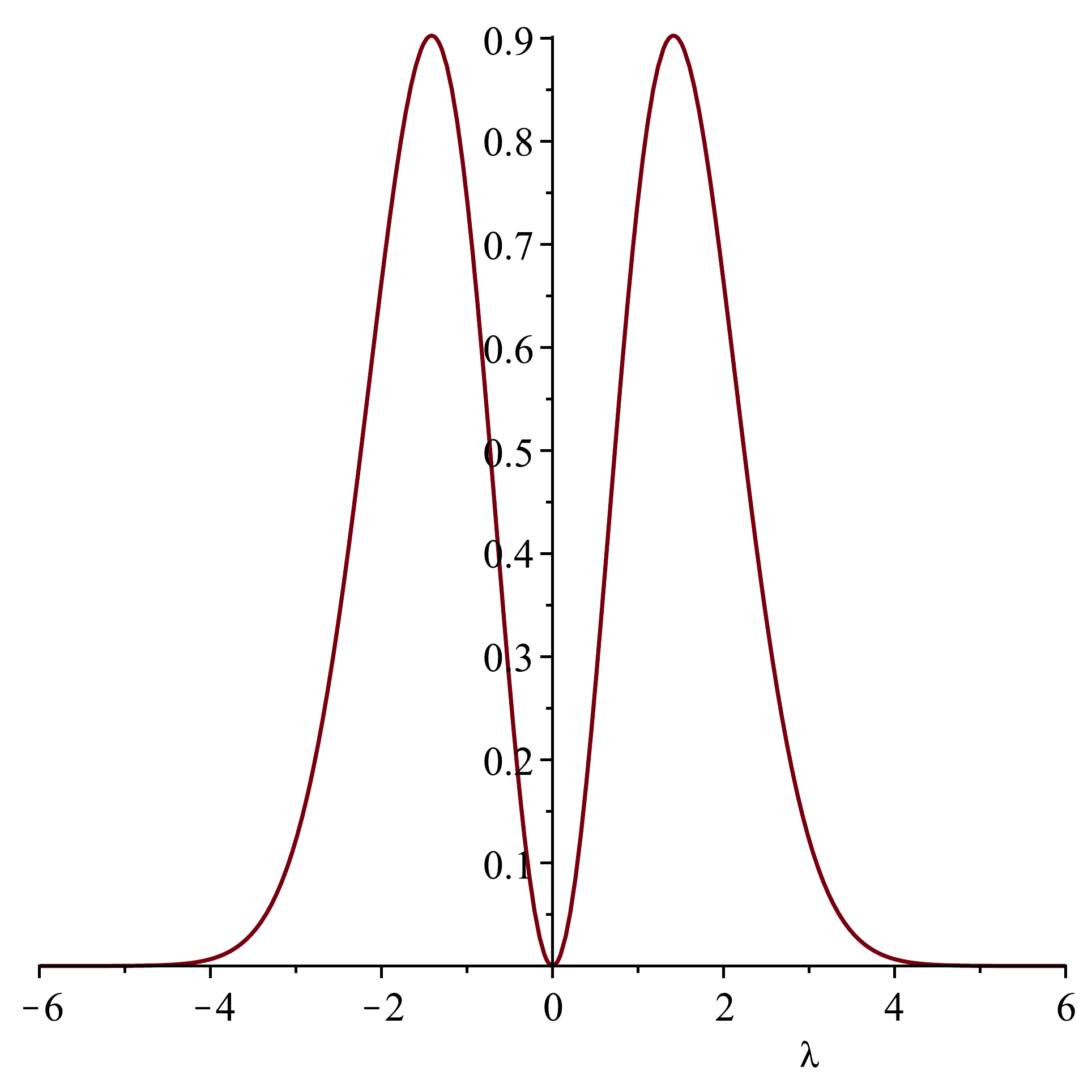}
 		\end{minipage}
 		       \caption{Plots of $\psi(t)$ and $\widehat \psi(\lambda).$ } \label{fig6}
 	\end{center}
 \end{figure}

The Mexican wavelet was used as a filter. It is defined by  \[ \psi(t)=\frac{2}{\sqrt{3\sigma}\pi^\frac{1}{4} }\left( 1-\left(\frac{t}{\sigma}\right)^2\right) \mathrm{e}^{-\frac{t^2}{2\sigma}}.\]    
Its Fourier transform  is, see \cite{Chun:2010}, \[\widehat \psi(\lambda) =\frac{\sqrt{8} \pi^\frac{1}{4} \sigma^\frac{5}{2}} {\sqrt{3}} \lambda^2 \mathrm{e}^{-\frac{\sigma^2 \lambda^2}{2}}.\]  Note, that the Fourier transform $\widehat \psi(\lambda)$  does not have a finite support but approaches zero very quickly when $\lambda\to +\infty.$ The value $\sigma=1$ was used in computations. Plots of $\psi (t)$ and $ \widehat \psi (\lambda) $ are shown in Figure \rm{\ref{fig6}}. In this case $c_2 $ and $c_3$ are 2 and 10, respectively.  

\begin{figure}[h]
 \begin{minipage}{0.5\textwidth}
         \centering
        \includegraphics[width= 7.5cm,height= 5cm,trim=1cm 1cm 1cm 0.5cm, clip=true] {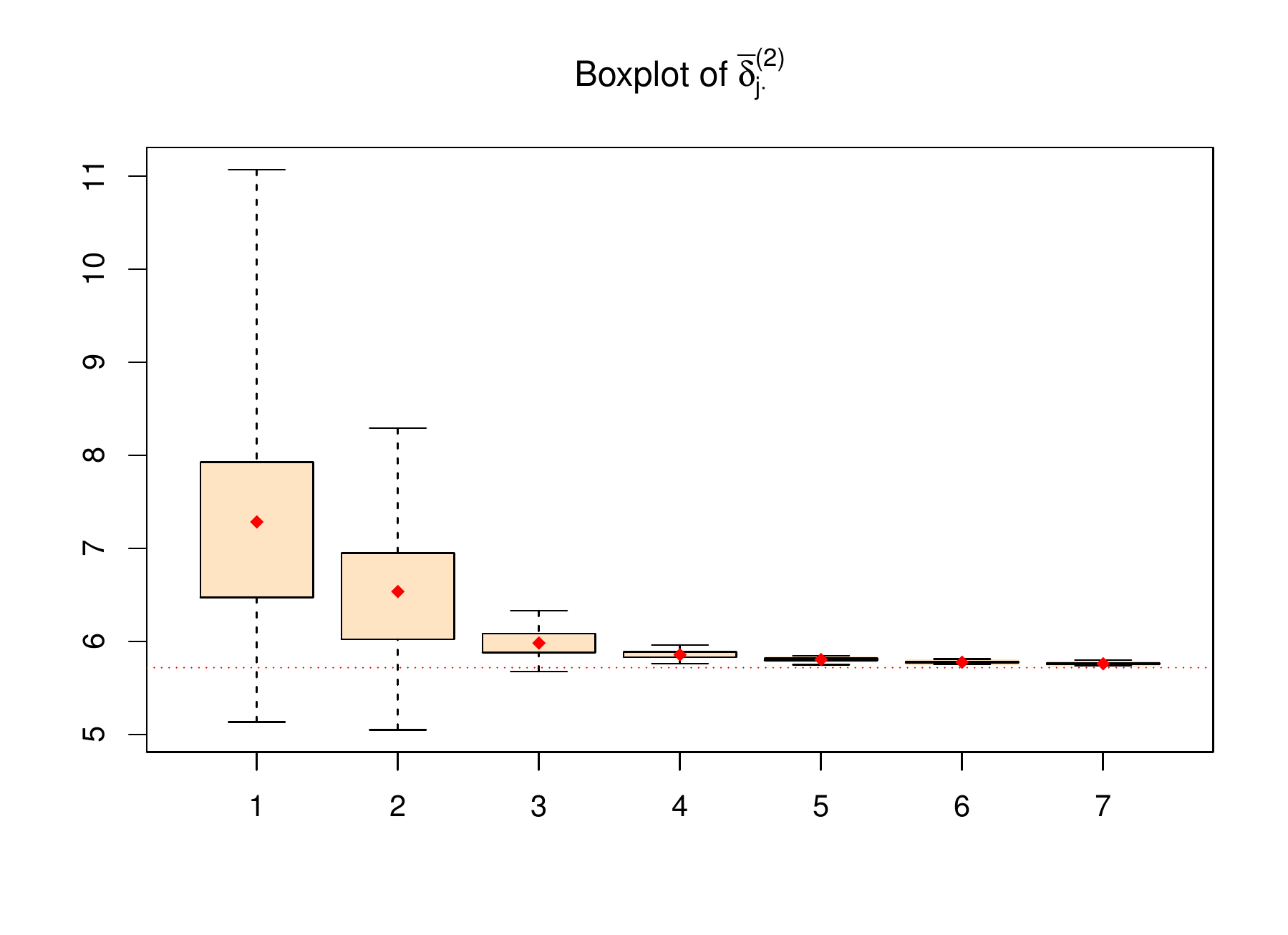}
\end{minipage}\hfill
\begin {minipage}{0.5\textwidth}
   \centering
  \includegraphics [width= 7.5cm,height= 5cm,trim=1cm 1cm 1cm 0.5cm, clip=true] {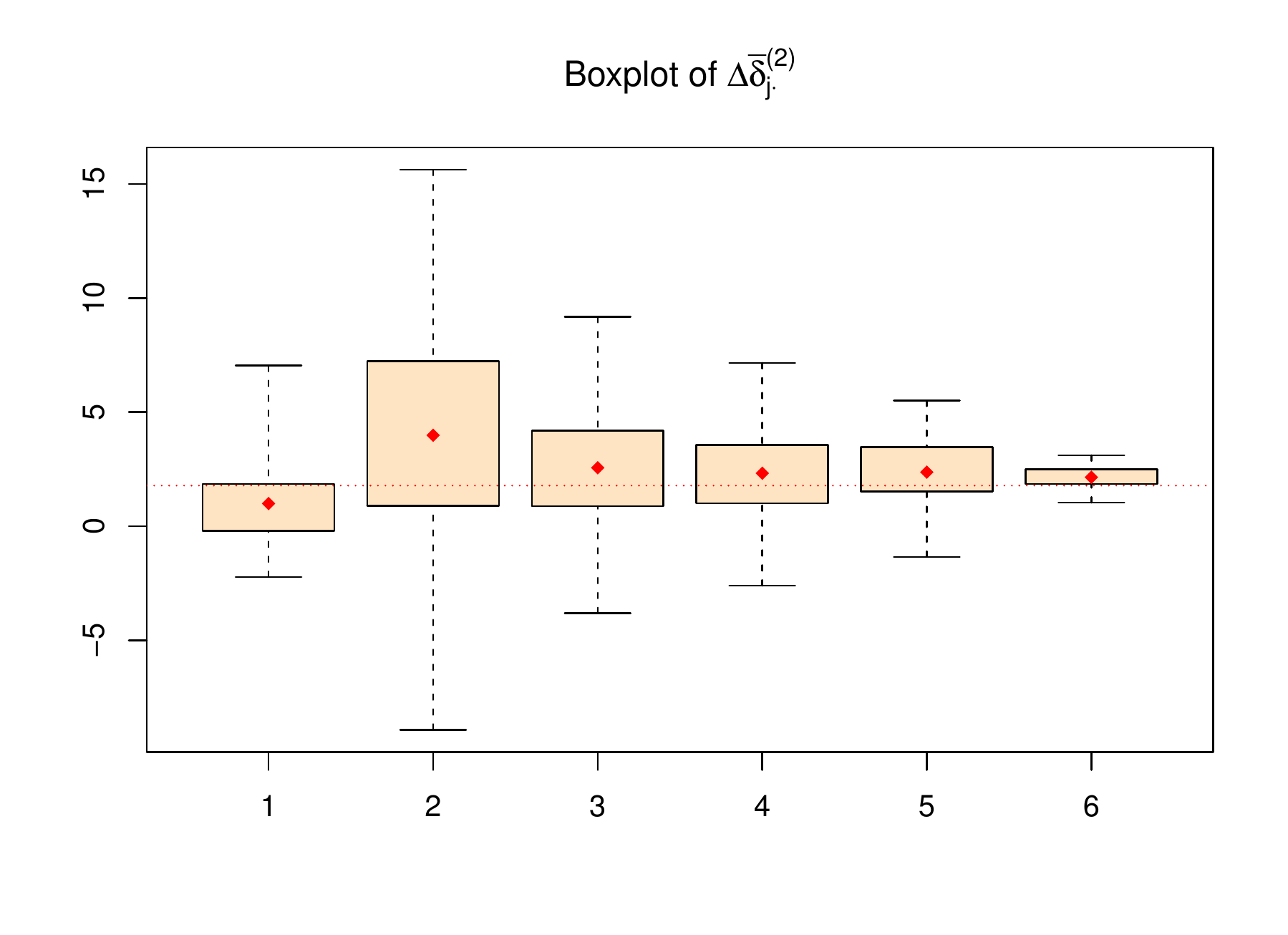}
\end{minipage}
        \caption{Boxplots of $\bar{\delta}_{j\cdot }^{(2)}$ and $\Delta \bar{\delta}_{ j\cdot}^{(2)}.$ } \label{boxplot:1}
\end{figure}

\begin{figure}[h]      
  \begin{minipage}{0.5\textwidth}
	\centering
	\includegraphics[width= 7.5cm,height= 5cm,trim=1cm 1cm 1cm 0.5cm, clip=true] {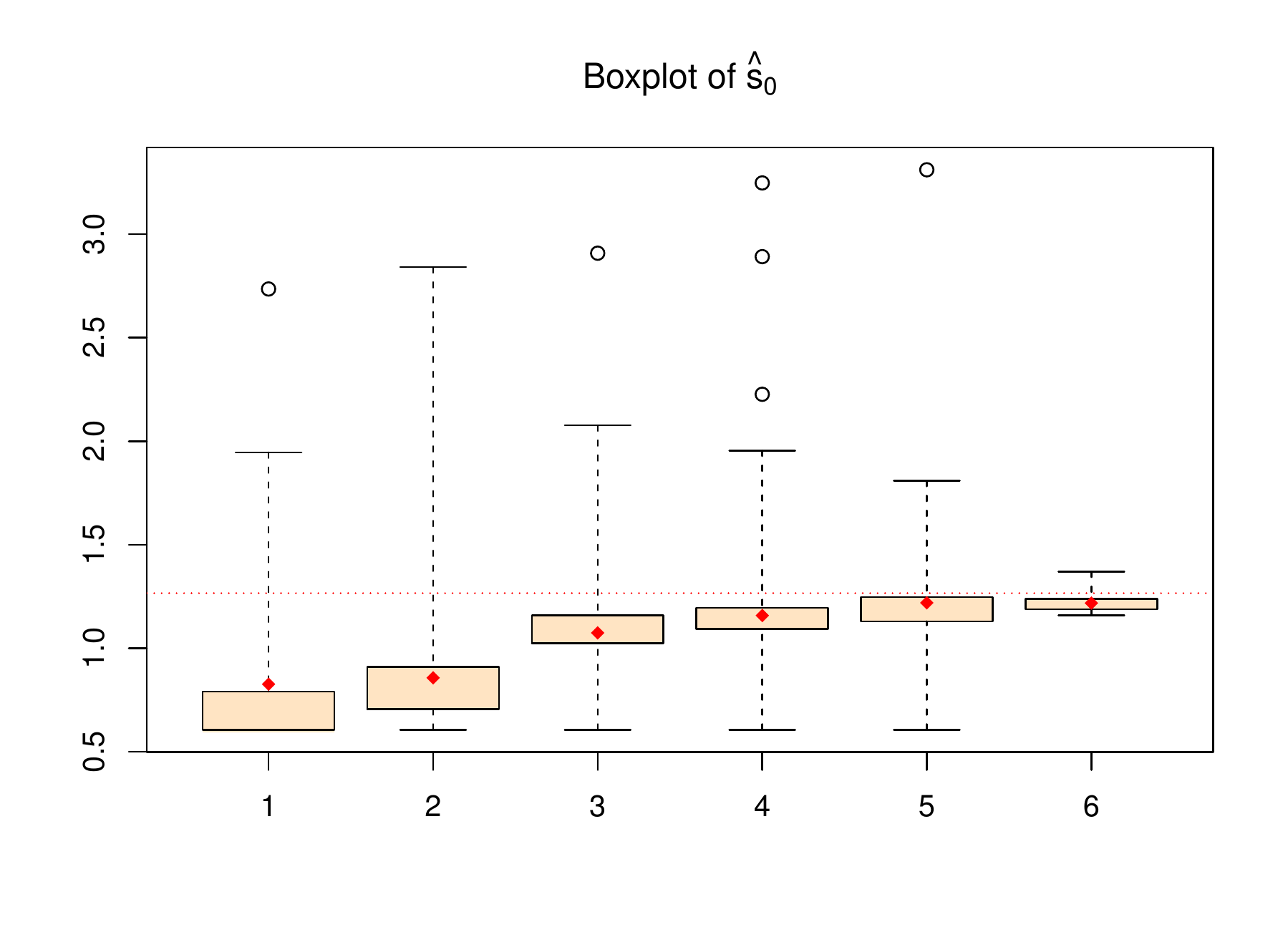}
\end{minipage}\hfill
\begin {minipage}{0.5\textwidth}
\centering
\includegraphics [width= 7.5cm,height= 5cm,trim=1cm 1cm 1cm 0.5cm, clip=true] {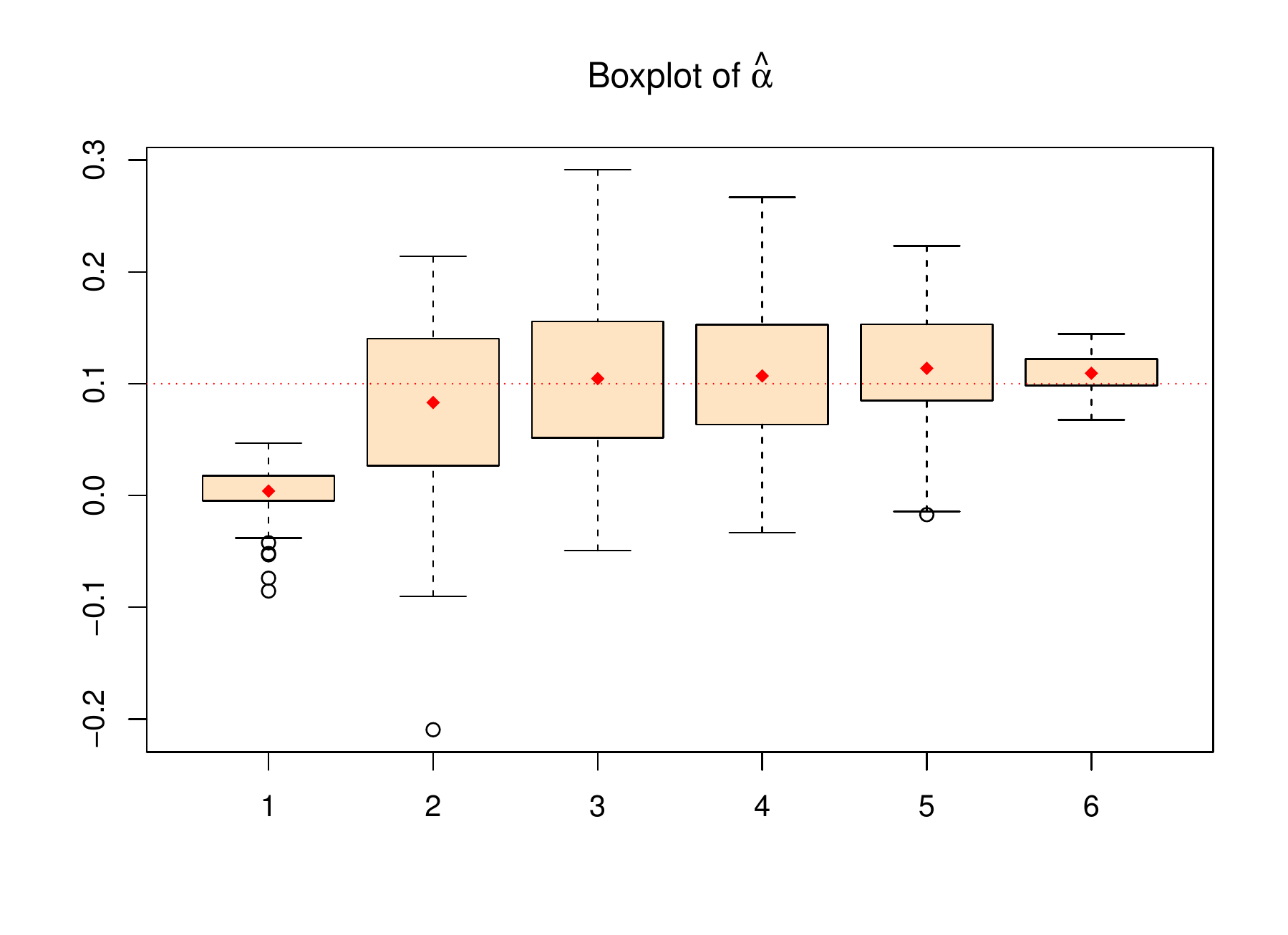}
\end{minipage}
       \caption{Boxplots of $\widehat s_{0j}$ and $\widehat \alpha_{j}.$} \label{boxplot:2}
\end{figure}

 100 realizations of $X(t)$ were generated and the corresponding wavelet coefficients $\delta_{jk}$ were calculated.  At first, to computed $\hat{s}_{0j}$ and $\hat{\alpha}_{j},$  the statistics  $ \bar{\delta}_{ j\cdot}^{(2)}$ and $\Delta \bar{\delta}_{ j\cdot}^{(2)}$ were found  using $a_j=j,$ $\gamma_j=1,$ $r_j=a_j^{-2.5}$ and $ m_j=a_j^9,$ $ j\geq1.$ These sequences satisfy the assumptions of Theorem~\ref{Theorem_1}.  Figure~\ref{boxplot:1} displays box plots of $\bar{\delta}_{ j\cdot}^{(2)}$ and $\Delta \bar{\delta}_{ j\cdot}^{(2)}$ for the simulated  realizations. The values of $j$ are shown along the horizontal axe. The horizontal dashed lines show the true values of the corresponding parameters. These plots confirm that  $\bar{\delta}_{j\cdot }^{(2)}$ and  $\Delta \bar{\delta}_{ j\cdot}^{(2)}$ converge as $j$ increases. As expected, consult the upper bound (\ref{4_3}) in Proposition~\ref{pro_2}, the rate of convergence of $\Delta \bar{\delta}_{ j\cdot}^{(2)}$ is slower than in the case of $\bar{\delta}_{j\cdot }^{(2)}$.  Finally, the estimates $\hat{s}_{0j}$ and $\hat{\alpha}_{j}$ were calculated by (\ref{5_6}) for each simulation. Figure~\ref{boxplot:2} demonstrates that $\hat{s}_{0j}$ convergence to $\ s_{0}=\arccos(u)$  and $\hat{\alpha}_{j}$ to $\alpha=d$ as $j$ increases. 
 
 As the true values  of  parameters correspond to a point inside of $R_\mathbf{y}$ the majority of the parameters estimates are in the admissible region. Figures~\ref{boxplot:1} and~\ref{boxplot:2}  suggest that the adjusted statistics should be applied mainly for the cases $j=1$ and 2.
 
The table below gives numerical values of mean squared errors (MSE) for each parameter estimated in Figures~\ref{boxplot:1} and~\ref{boxplot:2}. These results numerically confirm the theoretical convergence.

\begin{center}
\begin{tabular}{|l|c|c|c|c|c|c|}
	\hline \rule[-1ex]{0pt}{2.5ex}
	\backslashbox{$MSEs$}{$j$  } & 1 & 2 & 3 & 4 & 5 & 6 \\ 
	\hline \rule[-1ex]{0pt}{2.5ex}
	$MSE(\bar{\delta}_{j\cdot }^{(2)})$ & 3.5277 & 1.0751 & 0.0880 & 0.0214 & 0.0079 & 0.0035  \\ 
	\hline \rule[-1ex]{0pt}{2.5ex}
	$MSE(\Delta \bar{\delta}_{ j\cdot}^{(2)})$ &  2.9610 & 26.6364 &  8.1881 &  4.1055 &  2.2725 &  0.3252  \\ 
	\hline \rule[-1ex]{0pt}{2.5ex}
	$MSE(\hat{s}_{0j})$ &  0.3757 & 0.2509 & 0.1349 & 0.1524 & 0.0784 & 0.0038  \\ 
	\hline \rule[-1ex]{0pt}{2.5ex}
	$MSE(\hat{\alpha}_{j})$ & 0.0098 & 0.0066 & 0.0067 & 0.0046 & 0.0027 & 0.0003  \\ 
	\hline 
\end{tabular} 
\end{center}

\section{Directions for future research}
This paper has discussed statistical inference for parameters of seasonal long-memory processes with a spectral singularity at a non-zero frequency. The results were derived for wide classes of models with  Gegenbauer-type spectral densities using very general filter transforms.
 
An important area for future explorations is obtaining similar results for the case of multiple singularities with the long-memory parameters varying across singularity locations, see the discussion on SCLM
 (Seasonal/Cyclical Long Memory) in \cite{ArtRob:1999}. For the case of multiple unknown parameters one can derive additional estimation equations similar to the ones in Section \ref{sec_{4}} using higher order differences of $\bar{\delta}_{ j\cdot}^{(2)}.$ 

As this paper studied the case of Gegenbauer-type spectral densities given in Assumption \ref{Assumption_1}, it would be interesting to apply the developed methodology to other seasonal long-memory models.

This paper  develops statistical inference for parameters using functional data. There are numerous application where $X(t)$ is observed  only on a discrete grid or at random  moments of a finite time interval. Also, seasonal long-memory processes are often determined by discrete-time fractional autoregressive integrated moving average FARIMA models. In such cases approximate formulas are used to compute filter transforms, see Section 3.2 in \cite{Bardet:2010}. We plan to investigate statistical properties of the corresponding "approximate" estimates using approaches similar to \cite{Ayaber:2011} and \cite{Bardet:2010}. 

Finally, it is important to extend the methodology to the multidimensional case of random fields, see the discussion in \cite {Espejo:2015}.

\section*{Acknowledgments}
Andriy Olenko is grateful  to Laboratoire d'Excellence, Centre Europ\'{e}en pour les Math\'{e}matiques, la Physique et leurs interactions (CEMPI, ANR-11-LABX-0007-01), Laboratoire de Math\'{e}matiques Paul Painlev\'{e}, France, for support and giving him the opportunity to pursue research at the Universit\'{e} Lille 1 for two months.

Andriy Olenko was partially supported under the Australian Research Council's Discovery Projects funding scheme (project number  DP160101366) and the La Trobe University DRP Grant in Mathematical and Computing Sciences.\\

\textbf {Supplementary Materials} The codes used for simulations and examples in this article are available in the folder "Research materials" from \url{https://sites. goggle .com/site/olenkoandriy/}.

\end{document}